\documentclass[11pt,a4paper]{article}
\RequirePackage[OT1]{fontenc}
\usepackage[english]{babel}
\usepackage{amsmath}
\usepackage{amssymb}
\usepackage{bbm}
\usepackage{bm}
\usepackage{mathrsfs}
\usepackage{dsfont}
\usepackage{graphicx,theorem}
\usepackage{color}

\allowdisplaybreaks  
\newcommand{\X}{\mathbb{X}}

\newcommand{\indic}{\mathds{1}}

\newcommand{\R}{\mathds{R}}
\newcommand{\N}{\mathds{N}}

\newcommand{\ep}{\epsilon}

\newcommand{\Ccr}{\mathscr{C}}

\newcommand{\ddr}{\mathrm{d}}

\newcommand{\edr}{\mathrm{e}}

\newcommand{\p}{\tilde{p}}
\newcommand{\pr}{\mathrm{P}}
\newcommand{\E}{\mathrm{E}}

\newtheorem{T} {Theorem}

\newtheorem{Prop}[T]{Proposition}

\newtheorem{rem}[T]{Remark}
\def\simind{\stackrel{\mbox{\scriptsize{ind}}}{\sim}}
\def\simiid{\stackrel{\mbox{\scriptsize{iid}}}{\sim}}

\newcommand{\comillas}[1]{``\,#1\,"}

\allowdisplaybreaks

\long\def\symbolfootnote[#1]#2{\begingroup\def\thefootnote{\hspace*{-1mm}\fnsymbol{footnote}}\footnote[#1]{#2}\endgroup}

\usepackage[bookmarks,bookmarksnumbered,colorlinks=true,pdfstartview=FitV, linkcolor=red,citecolor=blue,urlcolor=blue]{hyperref}

\begin{document}

\begin{center}
\textmd{\Large{\bfseries{{Are Gibbs--type priors the most natural generalization of the Dirichlet process?\symbolfootnote[0]{(c) 2015 IEEE. Personal use of this material is permitted. Permission from IEEE must be obtained for all other users, including reprinting/republishing this material for advertising or promotional purposes, creating new collective works for resale or redistribution to servers or lists, or reuse of any copyrighted components of this work in other works. Publisher version \href{http://dx.doi.org/10.1109/TPAMI.2013.217}{DOI: 10.1109/TPAMI.2013.217}}}}}}
\end{center}
\begin{center}
{P. De Blasi$^1$, S. Favaro$^1$, A. Lijoi$^2$, R.H. Mena$^3$, I. Pr\"unster$^1$ and M. Ruggiero$^1$}
\end{center}

\begin{quote}
\begin{small}
\begin{center}
\noindent $^1$ Universit\`a degli Studi di Torino and Collegio Carlo Alberto, Torino, Italy.\\ \textit{E-mail}: pierpaolo.deblasi@unito.it; stefano.favaro@unito.it; igor.pruenster@unito.it; matteo.ruggiero@unito.it

\noindent $^2$ Universit\`a degli Studi di Pavia and Collegio Carlo Alberto, Torino, Italy.\\ \textit{E-mail}: lijoi@unipv.it

\noindent $^3$ Universidad Aut\'onoma de M\'exico, M\'exico\\ \textit{E-mail}: ramses@sigma.iimas.unam.mx

\end{center}
\end{small}
\end{quote}

\smallskip


\begin{abstract}

Discrete random probability measures and the exchangeable random
partitions they induce are key tools for addressing a variety of
estimation and prediction problems in Bayesian inference. Indeed, many
popular nonparametric priors, such as the Dirichlet and the
Pitman--Yor process priors, select discrete probability distributions
almost surely and, therefore, automatically induce exchangeable random
partitions. Here we focus on the family of \textit{Gibbs--type
  priors}, a recent and elegant generalization of the Dirichlet and the
Pitman--Yor process priors. These random probability measures share
 properties that are appealing both from a theoretical and an applied point
 of view: (i) they admit an intuitive characterization in terms of
their predictive structure justifying their use in terms of a precise
assumption on the learning mechanism; (ii) they stand out in terms of
mathematical tractability; (iii) they include several interesting special cases
besides the Dirichlet and the Pitman--Yor processes. The goal of our
paper is to provide a systematic and unified treatment of Gibbs--type
priors and highlight their implications 
for Bayesian
nonparametric inference. We will deal with their distributional
properties, the resulting estimators, frequentist asymptotic validation and the construction of time--dependent versions. Applications, mainly concerning hierarchical mixture models and species sampling, will serve to convey the main ideas. The intuition inherent to this class of priors and the neat results that can be deduced for it lead one to wonder whether it actually represents the most natural generalization of the Dirichlet process.
\par

\vspace{9pt}
\par
\end{abstract}


\section{Introduction and preliminaries}\label{sec1}

One of the main research lines within Bayesian Nonparametrics has been the proposal and study of classes of random probability measures whose laws act as nonparametric priors. Several such classes contain, as a special case, Ferguson's  Dirichlet process \cite{Fer(73)}, which  still represents the cornerstone of the field. A recent review that covers many of these models and uses completely random measures as a unifying concept can be found in \cite{lp}. When going beyond the Dirichlet
  process one typically has to face a trade--off between the desire of
  generality (which, as far as inference is concerned, implies
  flexibility of the model) and tractability, both analytical and
  computational. Probably the most successful proposal is represented
  by the two--parameter Poisson--Dirichlet process introduced in
  \cite{Per(92)} and further investigated in countless papers, most
  notably in \cite{P95,py}. See \cite{Pit(06)} for a comprehensive
  review from a probabilistic perspective. Such a process is also
  known as Pitman--Yor (PY) process, especially in the Machine Learning
  community, according to a terminology introduced in \cite{Ish01}
  which we will also adopt in the present paper. For our purposes it
  is important to note that the PY process reduces to the
  Dirichlet process by setting one of its parameters equal to
  $0$. Nonetheless, some
  important distributional features of the PY process are fundamentally different
  according as to whether the value of such a parameter is equal to $0$ or not.
  A clear understanding of this aspect is possible by identifying a large class of priors, which embeds the PY process as a
  special case. Such a class is given by Gibbs--type priors, introduced in \cite{Gne(05)} and only briefly addressed in the above mentioned review of nonparametric priors \cite{lp}, thus motivating the main focus of this paper.
In fact, by close inspection of the predictive
  structure they lead to, it will become apparent that the variety of
  distributional characteristics can be
  actually traced back to crucially different assumptions on the
  learning mechanism. This leads to a novel classification of discrete nonparametric priors which also serves as motivation for the use of Gibbs--type priors. Moreover, Gibbs--type priors have the advantage
  of pinning down, in a neat way, the analytic tractability
  issue related to
  general classes of nonparametric priors: in fact,
they allow to
  split the prediction rule in two stages and to highlight the key
  quantity allowing simplification of the relevant
  expressions. Indeed, throughout the following sections one can
  appreciate the beauty and simplicity of various analytical results
  that admit straightforward application to statistical
  inference. Finally, it is to be
  noted that Gibbs--type priors include other notable special cases of
  priors beyond the Dirichlet and the PY processes:
  for example, normalized inverse Gaussian processes \cite{Lij(05)}
  and their generalization given by normalized generalized gamma
  processes \cite{LMP07} as well as mixtures of symmetric Dirichlet distributions \cite{Gne(05)}.
  Given this, can one state with confidence that Gibbs--type priors are a natural generalization of the Dirichlet process, maybe the most natural?

The present paper aims at providing a survey on Gibbs--type priors that accounts for recent findings both in the probabilistic and statistical literature. This will serve as an important opportunity for pointing out their analytical tractability, flexibility and suitability in a variety of inferential problems beyond current applications which include mixture models (see, e.g., \cite{Ish01,LMP07}), linguistics and information retrieval in document modeling (\cite{Teh06,Teh10}), species sampling (\cite{Lij(07),Lij(07Bio),Nav08}) and survival analysis \cite{jar(10)}, among others.

\subsection{Discrete random probability measures, exchangeable random partitions and predictive distributions} \label{sec:eppf}

We first lay out the basics of Bayesian inference in an exchangeable
framework and focus on some key concepts and tools. Suppose $(X_n)_{n
  \ge 1}$ is an (ideally) infinite sequence of observations, with each
$X_i$ taking values in some set $\X$. Moreover, $\mathbf{P}_\X$ is the
set of all probability measures on $\X$. Assuming $(X_n)_{n \ge
  1}$ to be \textit{exchangeable} is equivalent to assuming the existence of a probability distribution $Q$ on $\mathbf{P}_\X$ such that
    \begin{equation}
      \begin{split}
        X_i\,|\,\tilde p \,
        &\simiid\, \tilde p, \qquad i=1,\ldots,n\\
        \tilde p \, &\sim\, Q
      \end{split}
      \label{eq:nonpmic}
    \end{equation}
for any $n\ge 1$. Hence, $\tilde p$ is a random probability measure on $\X$ and its probability distribution $Q$, also termed \textit{de~Finetti measure}, represents the prior distribution when \eqref{eq:nonpmic} is used as a Bayesian model with an observed sample $X_i$, $i=1,\ldots,n$. Whenever $Q$ degenerates on a finite dimensional subspace of $\mathbf{P}_\X$, the inferential problem is usually called \textit{parametric}. On the other hand, when the support of $Q$ is infinite--dimensional then one typically speaks of a \textit{nonparametric} inferential problem and it is generally agreed (see, e.g., \cite{Feg74}) that having a large topological support is a desirable property for a nonparametric prior.
Given a sample $X_1, \ldots, X_n$ generated through \eqref{eq:nonpmic}, the (one--step ahead) predictive distribution coincides with the posterior expected value of $\tilde p$, that is
    \begin{equation}\label{eq:pred_general}
    \pr(X_{n+1}\in \, \cdot \, | X_1, \ldots, X_n)=
    \int_{\mathbf{P}_\X} p(\, \cdot \,)\, Q(\ddr p\,|\, X_1, \ldots, X_n),
    \end{equation}
where $Q(\,\cdot\,|\,X_1, \ldots, X_n)$ denotes the posterior distribution of $\tilde p$.

Discrete nonparametric priors, i.e.~priors which select discrete
distributions with probability $1$, play a key role in most Bayesian
nonparametric procedures. It is well--known that the Dirichlet
and the PY process priors share this property and the same
can be said for the broader class of Gibbs--type priors.
In fact, any random probability measure associated to a discrete prior can be represented as
    \begin{equation}
      \label{eq:discrete}
     \tilde p=\sum_{j=1}^\infty \tilde p_j\:\delta_{Z_j}
    \end{equation}
where $\delta_c$ stands for the unit point mass concentrated at $c$, $(\tilde p_j)_{j\ge 1}$ is a sequence of non--negative random variables such that $\sum_{j\ge 1}\tilde p_j=1$, almost surely, and $(Z_j)_{j\ge 1}$ is a sequence of $\X$--valued random variables. Henceforth, we further assume that $(\tilde p_j)_{j\ge 1}$ and $(Z_j)_{j\ge 1}$ are independent and that the $Z_j$'s are iid from a diffuse probability measure $P^*$ on $\X$  (or in other terms $\pr(Z_i\ne Z_j)=1$ for any $i\ne j$). Such a general subclass of discrete random probability measures has been called \textit{species sampling models} by Pitman \cite{P95}, a terminology that will be clarified in the following section.

As far as the observables $X_i$'s are concerned, the discrete nature of $Q$ implies that any sample $X_1, \ldots, X_n$ will feature ties with positive probability, therefore generating $K_n=k\leq n$ distinct observations $X_1^*, \ldots, X_{k}^*$ with frequencies $n_1, \ldots, n_k$ such that $\sum_{i=1}^k n_i=n$. When choosing and analyzing specific predictive structures, the key quantity to consider, from both a conceptual and a mathematical point of view, is the probability of observing a new distinct value not included in the sample $X_1,\ldots,X_n$, namely
    \begin{equation}
    \pr(X_{n+1}=\hbox{``new''} \: |\: X_1, \ldots, X_n),
    \label{eq:pred_new}
    \end{equation}
which will appear throughout the paper.
To 
be more concrete consider the Dirichlet and the
PY processes. In the Dirichlet case, with parameters given by
$P^*$ and $\theta>0$, one has
    \begin{equation*}
    \pr(X_{n+1}=\hbox{``new''} \: |\: X_1, \ldots, X_n)=\frac{\theta}{\theta+n}.
    \label{eq:pred_new_dir}
    \end{equation*}
In the PY case, in addition to $P^*$, one has two parameters $(\sigma,\theta)$ whose admissible values are  $\sigma\in[0,1)$ with $\theta>-\sigma$ or $\sigma<0$ with $\theta=m|\sigma|$ for some positive integer $m$. One then has
    \begin{equation*}
    \pr(X_{n+1}=\hbox{``new''} \: |\: X_1, \ldots, X_n)=\frac{\theta+\sigma k}{\theta+n}
    \label{eq:pred_new_PY}
    \end{equation*}
from which it is apparent that the corresponding probability for the Dirichlet process is recovered by setting $\sigma=0$.

Within such a framework, discrete random probability measures can be characterized in terms of the exchangeable random partition they imply, another key aspect of the paper for which we provide some essential background. Given the discreteness of $Q$, $\tilde p$ induces a partition
of $X_1,\ldots,X_n$ that is well described by means of an extremely useful tool, namely the \textit{exchangeable partition probability function} (EPPF) \cite{P95} given by
    \begin{equation}
    \label{eq:eppf_start}
    p_k^{(n)}(n_1,\ldots,n_k)=\int_{\X^k}\E\left(\tilde p^{n_1}(\ddr x_1)\,\cdots\,\tilde p^{n_k}(\ddr x_k)\right).
    \end{equation}
It is also of simple interpretability: it essentially corresponds to the probability, induced by $\tilde p$, of observing a sample of size $n$, $X_1, \ldots, X_n$, exhibiting $K_n=k$ distinct observations with frequencies $n_1,\ldots,n_k$ or, equivalently, a specific partition into $K_n=k$ clusters with frequencies $n_1,\ldots,n_k$. See \cite{P95,Pit(06)} for details. Note also that an EPPF satisfies the addition rule
    \begin{equation}
    \label{eq:addition_rule}
    p_k^{(n)}(n_1,\ldots,n_k)=p_{k+1}^{(n+1)}(n_1,\ldots,n_k,1)+\sum_{j=1}^k p_k^{(n+1)}(n_1,\ldots,n_j+1,\ldots,n_k).
    \end{equation}
For both Dirichlet and PY processes, the EPPF is available in closed form. In the former case it is given by
\begin{equation}
  \label{eq:eppf_dirichlet}
  p_k^{(n)}(n_1,\ldots,n_k)=\frac{\theta^k}{(\theta)_n}\prod_{i=1}^k (n_i-1)!
\end{equation}
where $(\theta)_n=\theta(\theta+1)\,\cdots\,(\theta+n-1)$ for any
$n\ge 1$. For the PY process, it coincides with
\begin{equation}
  \label{eq:eppf_py}
  p_k^{(n)}(n_1,\ldots,n_k)=\frac{\prod_{i=1}^{k-1}(\theta+i\sigma)}{(\theta+1)_{n-1}}\:
\prod_{i=1}^k (1-\sigma)_{n_i-1}.
\end{equation}
The identification of the EPPF leads to the direct determination of the predictive distribution in \eqref{eq:pred_general}. Indeed, if $X_1,\ldots,X_n$ is a sample featuring $k\le n$ distinct values with respective frequencies $n_1,\ldots,n_k$, one has
\begin{equation}\label{eq:ss_new}
\pr(X_{n+1}=\hbox{``new''} \: |\: X_1, \ldots,
X_n) = \frac{p_{k+1}^{(n+1)}(n_1,\ldots,n_k,1)}{p_{k}^{(n)}(n_1,\ldots,n_k)}
\end{equation}
and the predictive distribution in \eqref{eq:pred_general} is a linear
combination of $P^*(\,\cdot\,)=\E(\p(\, \cdot \, ))$, which
can be interpreted as the prior guess at the shape of $\p$, and of a weighted measure of the observations, namely
\begin{multline}
  \label{eq:pred_eppf}
  \pr(X_{n+1}\in\,\cdot\,|X_1,\ldots,X_n)=\\
  \frac{p_{k+1}^{(n+1)}(n_1,\ldots,n_k,1)}{p_{k}^{(n)}(n_1,\ldots,n_k)}\:P^*(\,\cdot\,)
  +\sum_{j=1}^k
  \frac{p_{k}^{(n+1)}(n_1,\ldots,n_j+1,\ldots,n_k)}{p_{k}^{(n)}(n_1,\ldots,n_k)}
  \delta_{X_j^*}(\,\cdot\,).
\end{multline}
Note that the right hand side of \eqref{eq:pred_eppf} is guaranteed to sum up to $1$ if evaluated over the whole space $\X$ by the addition rule \eqref{eq:addition_rule}. In the PY process case the predictive distribution takes on the form
\begin{equation}
  \label{eq:pred_PY}
  \pr(X_{n+1}\in\,\cdot\,|X_1,\ldots,X_n)=\frac{\theta+\sigma k}{\theta+n}\:P^*(\,\cdot\,) +\frac{1}{\theta+n}\sum_{j=1}^k
  (n_j-\sigma) \delta_{X_j^*}(\,\cdot\,)
\end{equation}
which, for $\sigma=0$, reduces to the well--known Dirichlet process predictive structure given by a linear combination of $P^*$ and the empirical measure.

\subsection{Applications to species sampling problems and mixture modeling}

Discrete nonparametric priors in general, and Gibbs--type priors in particular, are suited for addressing inferential issues that arise in species sampling problems and in mixture modeling, among others. We now briefly sketch these frameworks. 

Consider a discrete random probability measure \eqref{eq:discrete} with the specifications as in Section \ref{sec:eppf}. It is, then, apparent that \eqref{eq:discrete} can be seen as a tool for describing the structure of a population made of different types or species with certain proportions, which are modeled through \eqref{eq:discrete} as random proportions $\tilde p_j$. On the basis of
this fact, one can equivalently use the $Z_i$'s or the positive integers $\{1,2,\ldots\}$ to label different species or types that can be sampled. Indeed, if $(\xi_n)_{n\ge 1}$ is an auxiliary integer--valued sequence such that $\pr(\xi_n=j\,|\,\tilde p)=\tilde p_j$, for any $n$ and $j$, model~\eqref{eq:nonpmic} corresponds to assuming that $X_i=Z_{\xi_i}$. Hence the $X_n$'s can be interpreted as the \textit{observed species labels} since, due to the diffuse nature of $P^*$, any two data points $X_i$ and $X_j$, for $i\ne j$, differ if and only if $\xi_i$ and $\xi_j$ do. Moreover, one has that $\pr(\xi_i=\xi_j)>0$, for any $i\ne j$, and this entails that the $i$--th and the $j$--th observations may reveal the same species with
positive probability. It is precisely this connection which motivates
the terminology adopted in \cite{Pit(96)}, \textit{species sampling
  model}. Moreover, an exchangeable sequence $(X_n)_{n\ge 1}$ for
which \eqref{eq:nonpmic} holds true, with $\tilde p$ a species sampling model, takes on the name of \textit{species sampling sequence}.

By virtue of this interpretation, there are a number
of statistical problems one can face adopting a Bayesian nonparametric
perspective. Indeed, in many statistical applications
one typically observes a sample of species labels
$X_1,\ldots,X_n$ and designs further sampling $X_{n+1},\ldots,X_{n+m}$
on the basis of estimates of some quantities of interest such as,
e.g.: the number of new distinct species that will be detected in a
new sample of size $m$; the number of species with a given frequency,
or with frequency below a certain threshold, in
$X_1,\ldots,X_{n+m}$; the probability that the $(n+m+1)$--th draw will
consist of a species having frequency $\ell\ge 0$ in
$X_1,\ldots,X_{n+m}$. These, in turn, provide measures of overall and
rare species diversity and are of interest in biological, ecological
or linguistic studies, just to mention a few. In this respect, the
predictive approach briefly sketched in Section~\ref{sec:eppf} plays
an important role and provides nice and elegant answers to these problems 
in the framework of Gibbs--type priors.

Discrete nonparametric priors are also basic building blocks for
hierarchical mixture models that are typically used for density
estimation and clustering but also in more complex dependent structures. To keep things simple consider the univariate
density estimation case and let $f(\,\cdot\,|\,\cdot\,)$ denote a
kernel defined on $\R\times\X$ and taking values in $\R^+$
such that $\int_\R f(y|x)\,\ddr y=1$, for any $x$ in $\X$. Hence,
$f(\,\cdot\,|x)$ defines a density function on $\R$, for any $x$. The
observations are then from a sequence $(Y_n)_{n\ge 1}$ of real--valued
random variables such that
\begin{equation}
  \label{eq:hierar}
  \begin{split}
    Y_i\,|\,X_i \: & \simind \: f(\,\cdot\,|X_i)\qquad i=1,\ldots,n\\
    X_i\,|\,\tilde p \: &\simiid\: \tilde p\qquad\qquad\quad i=1,\ldots,n\\
    \tilde p \: &\sim\: Q.
  \end{split}
\end{equation}
The typical choice for $\tilde p$ is represented by the Dirichlet
process leading to the Dirichlet process mixture model introduced by Lo \cite{lo}, which represents the most popular Bayesian nonparametric model to date. In addition to density estimation such model serves also clustering purposes. In fact, here $(X_n)_{n\ge 1}$ is a sequence of latent exchangeable random elements and the unobserved number $K_n$ of distinct values among $X_1,\ldots,X_n$ is the number of clusters into which the observations $Y_1,\ldots,Y_n$ can be grouped. Posterior inferences for $K_n$ are of great importance and the specification of a Gibbs--type prior $\tilde p$ in \eqref{eq:hierar} allows for an effective detection of the number of clusters that have generated the data.

\subsection{Outline of the paper}

Section~\ref{sec:ramses} first provides an intuitive characterization of Gibbs--type priors based on a suitable classification of species sampling models. This is, then, followed by a formal definition and an overview of their distributional properties that are of interest for applications
to Bayesian inference. Particular emphasis is given to the role
played by one of the parameters that characterizes them. Section~\ref{sec:mixtures}
discusses the use of Gibbs--type priors within hierarchical mixture
models for density estimation and clustering. Section~\ref{stefano} focuses on the application of Gibbs--type priors to prediction problems and Section~\ref{sec:consistency} deals with their frequentist asymptotic properties. Section~\ref{sec:matteo} concisely discusses extensions of Gibbs--type priors to dynamic contexts. Finally, Section~\ref{sec:conclusions} contains some concluding remarks trying to answer the question posed in the title of the paper.


\bigskip

\section{Gibbs--type priors}\label{sec:ramses}

An interesting and useful classification of species sampling models can be given in terms of the structure of the 
probability of generating a new value they induce. This leads to an intuitive characterization of Gibbs--type priors and represents also one of the main motivations for their use. Our result is somehow in the spirit of Zabell's \cite{zab(82)} characterization of the Dirichlet process in terms of the so--called \textit{Johnson's sufficientness postulate}.
To this end, recall that the key quantity is \eqref{eq:pred_new} representing the probability of generating a new value given the past associated to a species sampling model as specified in Section \ref{sec:eppf}. According to its structure one can classify the underlying priors in three main categories. Denote by $\bm{\Theta}$ a finite--dimensional parameter possibly entering the specification of $\p$ in \eqref{eq:nonpmic}. In general, one has $\pr(X_{n+1}=\hbox{``new''} \,|\, X_1, \ldots, X_n)=f(n, k, n_1, \ldots, n_k, \bm{\Theta})$, which means that the probability of obtaining a new observation depends on the sample size $n$, the number of distinct values $k$, their frequencies $(n_1, \ldots, n_k)$ and the parameter $\Theta$. We will denote $\pr(X_{n+1}=\hbox{``new''} \,|\, X_1, \ldots, X_n)$ by $f(n, k, \bm{\Theta})$  if it does not depend on $(n_1, \ldots, n_k)$ and by $f(n, \bm{\Theta})$ if it depends neither on $(n_1, \ldots, n_k)$ nor on $k$.

\medskip

\begin{Prop} \label{prop:classification} Let $\tilde p$ be a species sampling model. Then the following classification in terms of the structure of the probability of generating a new value holds:
\begin{itemize}
\item[(i)] $\pr(X_{n+1}=\hbox{``new''} \,|\, X_1, \ldots, X_n)=f(n, \bm{\Theta})$ if and only if $\tilde p$ is a Dirichlet process;
\item[(ii)] $\pr(X_{n+1}=\hbox{``new''} \,|\, X_1, \ldots, X_n)=f(n,
  k, \bm{\Theta})$ if and only if $\tilde p$ is of Gibbs--type;
\item[(iii)] $\pr(X_{n+1}=\hbox{``new''} \,|\, X_1, \ldots, X_n)=f(n, k, n_1, \ldots, n_k, \bm{\Theta})$ otherwise.
\end{itemize}
\end{Prop}

\medskip

Even if the Dirichlet process has proven to perform well in several applied contexts,
from a merely conceptual point of view it seems
  too restrictive to let the probability of generating new values
  depend solely on the sample size $n$ and on its total mass parameter $\theta$ and to summarize
  all other information contained in the data by means of a suitable
  specification of the scalar
  parameter $\theta$. One would like indeed such a probability to explicitly depend also on (at least) the number of distinct observed values, since it summarizes the heterogeneity in the sample. By virtue of (ii), this is tantamount to resorting to a Gibbs--type prior. According to the specific situation, one might
want to model \eqref{eq:pred_new} as an increasing or decreasing
function of $K_n$, which will be shown to correspond to
Gibbs--type priors with a parameter, to be identified later, being either positive or  negative, respectively. Case (iii), which corresponds to the most general setup and prediction of new values explicitly depends on all the information conveyed by the data, is in principle the most desirable prediction structure. However, there are two main operational problems that one needs to take into account. On the one hand, the general case (iii) gives rise to serious analytical hurdles
and priors have to be studied on a case-by-case basis typically leading to quite
complicated expressions (see \cite{Fav(11)}). On the other hand, it is not clear how one should explicitly specify the dependence of the probability of observing a new species on the observed frequencies $n_1,\ldots,n_k$ so that it reflects an opinion on the learning mechanism for the data. It is thus reasonable that such prior opinion be encoded through the finite--dimensional parameter $\bm{\Theta}$.
Hence, the above classification neatly shows the origin of the mathematical
tractability of Gibbs--type priors, which is due to a precise
simplifying assumption on the prediction structure. Overall, such an
assumption appears to be a satisfactory compromise between generality
(or flexibility) and tractability, and therefore motivates the attempt to study and understand the behavior of such priors.

After having stated and discussed a predictive characterization of
Gibbs--type priors, we now provide a different, though equivalent,
definition which is more useful when one wishes to analyze their distributional properties. As seen in Section~\ref{sec1}, a discrete nonparametric prior $\tilde p$  associated to an exchangeable scheme of the type \eqref{eq:nonpmic} can be characterized in terms of the associated EPPF $\{p_k^{(n)}:\: n\ge 1,\: 1,\le k\le n\}$ defined as in \eqref{eq:eppf_start}. Accordingly, one defines a \emph{Gibbs--type  prior} as a species sampling model such that
\begin{equation}
  \label{eq:1}
  p_k^{(n)}(n_1,\ldots,n_k)=V_{n,k}\:\prod_{i=1}^k (1-\sigma)_{n_i-1}
\end{equation}
for any $n\ge 1$, $k\le n$ and positive integers $n_1,\ldots,n_k$ such that $\sum_{i=1}^k n_i=n$, where $\sigma<1$ and the set of non--negative weights $\{V_{n,k}: \; n\ge 1, \; 1\le k\le n\}$ satisfies the forward recursive equation
    \begin{equation} \label{eq:recursion}
    V_{n,k}=(n-\sigma k) V_{n+1,k} +V_{n+1,k+1}
    \end{equation}
for any $k=1, \ldots, n$ and $n\ge 1$, with $V_{1,1}=1$. In light of
\eqref{eq:1} one can rephrase the reason for their tractability in
more mathematical terms, namely the product form of their
EPPFs which allows to handle conveniently the frequencies $n_i$. Given
\eqref{eq:1}, the probability of obtaining a new distinct observation
conditional on a sample $X_1,\ldots,X_n$ such that
  $K_n=k$ is
    \[
    \pr(X_{n+1}=\comillas{\mbox{new}}\mid X_1,\ldots,X_n)=\frac{V_{n+1,k+1}}{V_{n,k}}=f(n,k,\bm{\Theta})
    \]
which is in accordance with the above characterization.\\

\begin{rem}
According to the classification implied by Proposition~\ref{prop:classification}, mixtures of the Dirichlet process, obtained by mixing with respect to the total mass $\theta$ of the base measure, are in class (ii). To see this, let $\pi$ denote the prior on $\theta$ so that
$    \pi(\ddr \theta|X_1, \ldots, X_n) \propto \theta^k \pi(\ddr \theta)/(\theta)_{n}$,
where $(\theta)_n$ is the $n$--th ascending factorial. Hence
\[
    \mathbb{P}(X_{n+1}=\hbox{``new''} \: |\: X_1, \ldots, X_n)=\int_{\mathbb R^+}\frac{\theta^{k+1}}{(\theta)_{n+1}}\, \pi(\ddr \theta)
\]
will now depend on $k$. More generally, mixtures of Gibbs--type priors obtained by mixing with respect to a possible parameter entering the definition of $V_{n,k}$ are still of Gibbs--type and, thus, still lie in (ii). In contrast, Gibbs--type priors mixed with respect to $\sigma$ are such that $\pi(\ddr \sigma|X_1, \ldots, X_n) \propto V_{n,k}\:\prod_{i=1}^k (1-\sigma)_{n_i-1}\: \pi(\ddr \sigma)$ for some prior $\pi$ on $\sigma$. This clearly implies that the resulting family of species sampling models is in (iii), although one still preserves a Gibbs structure conditionally on $\sigma$.\\
\end{rem}

The definition \eqref{eq:1} implies that the induced predictive distributions are
    \begin{equation}
    \label{eq:predict_gibbs}
    \pr\left(X_{n+1}\in \cdot\: \big|\: X_1,\ldots,X_n \right)=
    \frac{V_{n+1,k+1}}{V_{n,k}}\, P^*(\cdot
    )+\frac{V_{n+1,k}}{V_{n,k}}\,\sum_{i=1}^k(n_i-\sigma)\,\delta_{X_i^*}(\cdot).
    \end{equation}
Hence, the predictive distribution is a linear convex combination of the prior guess $P^*$ at the shape of $\tilde p$ and of the weighted empirical distribution
$\widehat{P}_n=(n-k\sigma)^{-1}\sum_{i=1}^k (n_i-\sigma)\,\delta_{X_i^*}$.
The predictive structure \eqref{eq:predict_gibbs} exhibits some
appealing and intuitive features. In particular, the mechanism for
allocating the predictive mass among ``new'' and previously observed
data can be split into two stages. Given a sample $X_1,\ldots,X_n$,
the first step consists in allocating the mass between a new value
$X_{k+1}^*$ sampled from $P^*$ and the set of observed values
$\{X_1^*,\ldots,X_k^*\}$. This first step depends only on $n$ and $k$
and not on the frequencies $n_1, \ldots, n_k$. The second step is the following: conditionally on $X_{n+1}$ being a new value, it is sampled from the base measure $P^*$, whereas if $X_{n+1}$ coincides with one of the previously observed values $X_i^*$, for $i=1, \ldots, k$,
the coincidence probabilities are determined by the size $n_i$ of each
cluster and by $\sigma$. Hence, even if the
  frequencies $n_i$ do not affect the probability of
  allocating a predicted value between ``new'' and ``old'', they are
  explicitly involved conditional on the predicted
  value coinciding with a previously observed one: the more often a
  past observation is detected, the higher the probability of re--observing it. Also $\sigma$ plays an interesting role in weighting the empirical measure since, for $\sigma>0$, a reinforcement mechanism driven by $\sigma$ takes place. Indeed, one can see that the ratio of the probabilities assigned to any pair of $(X_i^*,X_j^*)$ is given by $(n_i-\sigma)/(n_j-\sigma)$. As $\sigma\to 0$, the previous quantity reduces to the ratio of the sizes of the two clusters and therefore
the coincidence probability is proportional to the size of the
cluster. On the other hand, if $\sigma>0$ and $n_i>n_j$, the ratio is
an increasing function of $\sigma$. Hence, as $\sigma$ increases the
mass is reallocated from $X_j^*$ to $X_i^*$. This means that the
sampling procedure tends to reinforce, among the observed clusters,
those having higher frequencies, which represents an appealing feature
in certain inferential contexts. See \cite{LMP07} for a discussion of
such reinforcement mechanisms and their use in mixture models. If
$\sigma<0$, the reinforcement mechanism works in the opposite way in the sense that the coincidence probabilities are less than proportional to the cluster size.

Besides influencing the balancedness of the partition of the exchangeable random elements directed by a Gibbs--type prior, the parameter $\sigma$ also determines the rate at which the number of clusters $K_n$ increases, as the sample size $n$ increases. As shown, e.g., in \cite{P03}, if we introduce
\begin{equation*}\label{K-n rescaling}
c_n(\sigma)=\left\{
  \begin{array}{ll}
    1 & \:\sigma<0\\
    \log n & \:\sigma=0\\
    n^\sigma & \:\sigma\in(0,1)
  \end{array}
  \right.
\end{equation*}
for any $n\ge 1$, then
\begin{equation}
  \label{eq:kappan_asympt}
  \frac{K_n}{c_n(\sigma)}\:\stackrel{\mbox{
      \begin{footnotesize}
        a.s.
      \end{footnotesize}}}{\longrightarrow}\: S_\sigma
\end{equation}
as $n\to\infty$. The limiting random variable $S_\sigma$ is termed
$\sigma$--\textit{diversity}. See \cite{Pit(06)} for details. It is
worth noting that if $\tilde p$ is the Dirichlet process with parameter measure $\theta P^*$, the
$\sigma$--diversity is degenerate on the total mass $\theta>0$ and
$K_n\sim \theta\,\log n$, for $n$ large enough, almost surely. This
special case was pointed out in \cite{KH73}. The larger $\sigma$, the faster the rate of increase of $K_n$ or, in other terms, the more new values are generated. Clearly, the case where $\sigma<0$ corresponds to a model accommodating for a finite number of distinct species in the population.

The combined effect of the reinforcement mechanism and the increase in
the rate at which new values are generated, both driven by $\sigma$,
is best visualized by looking at the special case of the PY
process. By close inspection of their predictive distributions
\eqref{eq:pred_PY} one notes that a new value, thus with frequency
$1$, entering the conditioning sample produces two
  effects: it is assigned a mass proportional to $(1-\sigma)$, instead of
  $1$, in the empirical component of the predictive and,
correspondingly, a mass proportional to $\sigma$ is added to the
probability of generating a new value. Therefore,
  if $\sigma>0$, new values are assigned a mass which is less than
  proportional to their cluster size (that is $1$) and the remaining
  mass is added to the probability of generating a new value. The
  first phenomenon gives rise to the reinforcement mechanism described
  above: if the new value is, then, re--observed it
    increases the associated mass by a quantity which is now proportional to
    $1$, and not less than proportional.
  The second effect implies that if $X_{n+1}$ is new, the probability
  of generating yet another new value, which overall still decreases
  as a function of $n$, is increased by a factor of
  $\sigma/(\theta+n+1)$. To sum up, the larger $\sigma$ the stronger
  is the reinforcement mechanism and at the same time the higher is the probability of generating a new value, which intuitively explains why one then obtains a growth rate of $n^\sigma$ for $K_n$. If $\sigma<0$ things work the other way round and one sees that each new generated value decreases the probability of generating further new values, thus providing intuition for the fact that in the end only a finite number of values will be generated. If $\sigma=0$, which corresponds to the Dirichlet process and mixtures of the Dirichlet process over the parameter $\theta$, everything is proportional to the cluster sizes which do not alter the probability of generating new values. As for another instance of a Gibbs--type prior, namely the normalized generalized gamma process that will be discussed later, a mechanism analogous to the PY process with $\sigma \in (0,1)$ can be identified though the proportionality constants that rescale the masses 
  are different due to the difference of the underlying $V_{n,k}$'s.

\subsection{Connections between Gibbs--type priors and product partition models}

There is also a close connection between Gibbs--type priors, and in particular the random partitions they induce, and exchangeable product partition models. The latter were introduced by \cite{Har(90)} and further studied, among others, by \cite{barry, quintana}. If $\Pi_n$ represents a random partition of the set of integers $\{1,\ldots,n\}$, a product
partition model corresponds to a probability distribution for $\Pi_n$ represented as follows
    \begin{equation}\label{eq:prod}
    \pr(\Pi_n=\{S_1,\ldots,S_k\})\propto \,\prod_{i=1}^k \rho(S_i)
    \end{equation}
where $\rho(\,\cdot\,)$ is termed \textit{cohesion function}. Now, let $|S|=\mbox{card}(S)$ and impose the cohesion function $\rho(\,\cdot\,)$ to depend only on the cardinality of the set $S$, that is $\rho(S_i):=\rho(|S_i|)=\rho(n_i)$. This is a natural and reasonable choice for a cohesion function. Then the random partition in \eqref{eq:prod} is, for any $n\geq k \geq 1$, the random partition induced by an exchangeable sequence if and only if $\rho(n_i)=(1-\sigma)_{n_i-1}/n_i!$ for $i=1, \ldots, k$ and $\sigma \in [-\infty, 1]$ with the proviso that $(1-\sigma)_{n_i-1}=1$ when $\sigma=-\infty$ and that $\Pi_n$ reduces to the singleton partition when $\sigma=1$. This is equivalent to saying that $\Pi_n$ is of Gibbs--type. Such a statement follows immediately from \cite{Gne(05)}. Therefore, random probability measures inducing exchangeable product partition models with cohesion function depending on the cardinality, i.e.
    \begin{align*}
    X_i^*|\Pi_n & \:\simiid\: P^* \qquad\qquad i=1,\ldots,K_n \\
    \Pi_n & \:\sim\: \mbox{product partition distribution with }
    \rho(S)=\rho(|S|),
    \end{align*}
coincide with the family of Gibbs--type priors.

\subsection{Sub--classes of Gibbs--type priors}\label{sec:sub}

Many nonparametric priors currently used for Bayesian inference represent particular cases of Gibbs--type priors, such as the Dirichlet process and the PY family. Indeed, it can be verified that
the set of weights
    \begin{equation}\label{vnkPD}
    V_{n,k}=\frac{\prod_{i=1}^{k-1}(\theta+i\sigma)}{(\theta+1)_{n-1}}
    \end{equation}
satisfies the recursive equation \eqref{eq:recursion} if the pair $(\sigma,\theta)$ is such that $\sigma\in[0,1)$ and $\theta>-\sigma$ or $\sigma<0$ and $\theta=m|\sigma|$ for some positive integer $m$. These constraints identify the set of admissible values of the parameters $(\sigma,\theta)$. The
corresponding Gibbs--type prior, identified by its EPPF \eqref{eq:1},
reduces to \eqref{eq:eppf_dirichlet} for  $\sigma=0$ therefore leading
to the Dirichlet process.
For any admissible $(\sigma,\theta)$ the EPPF \eqref{eq:1} coincides with \eqref{eq:eppf_py}, thus recovering the PY family. Another interesting special case of the PY process, and \textit{a fortiori} of Gibbs--type priors, is represented by the normalized  $\sigma$-stable process introduced by \cite{K75} which is obtained as a PY process with $\sigma\in(0,1)$ and $\theta=0$.

Before discussing other special cases of Gibbs--type priors, it is worth having a closer look at the PY family with $\sigma<0$ and $\theta=m|\sigma|$. In this case the weights in \eqref{vnkPD} are as follows
    \begin{equation}\label{vnkPD_1}
    V_{n,k}=\frac{|\sigma|^{k-1} \prod_{i=1}^{k-1}(m-i)}{(m|\sigma|+1)_{n-1}}\:\indic_{\{1,\ldots,\min(n,m)\}}(k)
    \end{equation}
where $\indic_A$ denotes the indicator function of set $A$. From \eqref{vnkPD_1} it is then easy to see (cfr. \cite{Pit(06)}) that the PY family with $\sigma<0$ and $\theta=m|\sigma|$ corresponds to a population composed by $m$ different species with proportions distributed according to a symmetric Dirichlet distribution with density function
    \[
    f_m(p_1,\ldots,p_{m-1})=\frac{\Gamma(m|\sigma|)}{\Gamma^m(|\sigma|)} \prod_{i=1}^{m-1}p_i^{|\sigma|-1}\:
    (1-p_1-\,\cdots\,p_{m-1})^{|\sigma|-1}
    \]
for any $(p_1,\ldots,p_{m-1})$ such that $p_i\ge 0$ for any $i$ and $\sum_{i=1}^{m-1}p_i\le 1$. Such a model arises, in the Population Genetics literature, as the stationary law of a Wright--Fisher model.

The PY family with parameters $(\sigma,m|\sigma|)$ and $\sigma<0$  is not only a distinguished
special case of Gibbs--type prior with $\sigma<0$ but actually is its
basic building block. In fact, any Gibbs--type random probability measure with
$\sigma<0$ is obtained by specifying a prior $\pi$ for the total number
of species $m$ in \eqref{vnkPD_1} and coincides with a species
sampling model having a random (finite) number of species. Crucially,
by \cite{Gne(05)}, the reverse implication holds true as
well: any Gibbs--type prior with $\sigma<0$ is a mixture of
PY processes with parameters $(\sigma, m|\sigma|)$, the
mixing measure being a probability measure on the positive
integers. Therefore, one can equivalently describe Gibbs--type priors with $\sigma<0$ in terms of a mixture model as
    \begin{equation}\label{eq:hier}
    \begin{array}{rcl}
    (\p_1,\ldots,\p_{\tilde m-1})\,|\, \tilde m &\sim&f_{\tilde m}\\
    \tilde m&\sim&\pi.
    \end{array}
    \end{equation}
Interesting special cases arise by particular specifications of $\pi$. For instance, if
    \begin{equation}\label{eq:mixing_gnedin}
    \pi(m)={\gamma(1-\gamma)_{m-1}\over m!}
    \end{equation}
for $m=1,2,\ldots$ with $\gamma\in(0,1)$, one obtains the model introduced by Gnedin \cite{Gne10}, which in the case of $\sigma=-1$ admits a completely explicit expression of the weights, namely
    \begin{equation}\label{eq:gnedin_model}
    V_{n,k}=\frac{(k-1)!(1-\gamma)_{k-1}(\gamma)_{n-k}}
    {(n-1)!(1+\gamma)_{n-1}}.
    \end{equation}
The peculiar feature of such a model, which makes it of great use in applications, is that the heavy-tailedness of \eqref{eq:mixing_gnedin} implies a model with finite random number of species whose expected value is infinite. Other interesting models are obtained by specifying the mixing distribution as a Poisson distribution restricted to the positive integers with parameter $\lambda>0$ , i.e.
  \begin{equation} \label{eq:mixing_poisson}
  \pi(m)={\edr^{-\lambda}\over 1-\edr^{-\lambda}} {\lambda^m\over m!}
  \end{equation}
for $m=1,2,\ldots$, or as a geometric mixing distribution
  \begin{equation}\label{eq:mixing_geometric}
  \pi(m)=(1-\eta)\eta^{m-1}
  \end{equation}
for some $\eta\in(0,1)$ and $m=1,2,\ldots$. These will be further discussed in Section~\ref{sec:cons_ex}.

Another important sub--class of Gibbs--type priors is the normalized generalized Gamma (NGG) process which corresponds to
  \begin{equation}\label{vnkGG}
    V_{n,k}=\frac{\mathrm{e}^{\beta}\: \sigma^{k-1}}{\Gamma(n)} \, \sum_{i=0}^{n-1}\binom{n-1}{i}\,(-1)^i\, \beta^{i/\sigma}\, \Gamma\left(k-\frac{i}{\sigma};\,\beta\right),
 \end{equation}
where $\sigma\in(0,1)$, $\beta>0$ and $\Gamma(x,a)=\int_x^\infty s^{a-1}\,\edr^{-s}\,\ddr s$ is the incomplete gamma function. Also the NGG process contains several interesting special cases: if $\sigma\to 0$ one obtains the Dirichlet process, whereas $\sigma=1/2$ yields the normalized inverse Gaussian process (N--IG) of \cite{Lij(05)}, which stands out for the availability of a closed form expression of its finite--dimensional distributions. Furthermore, if $\beta=0$, the normalized $\sigma$-stable process is also recovered from the NGG. See \cite{rlp,james2, LMP07}. The name attributed to this particular Gibbs--type prior is motivated by the fact that it can be defined by normalizing a generalized gamma completely random measure introduced in \cite{brix} and it therefore also belongs to the class of normalized random measures with independent increments (NRMI) introduced in \cite{rlp}. Interestingly, as shown in \cite{lpw2}, it turns out to be the only random probability measure belonging to both classes, NRMIs and Gibbs--type priors. All other NRMIs, such as for instance the generalized Dirichlet process in \cite{Lij(05b),Fav(11)}, are not of Gibbs--type.

In addition to specific examples described so far and still for the case of $\sigma>0$, one might wonder whether starting from the prediction rules \eqref{eq:predict_gibbs} it is possible to identify the class of random probability measures generating them. The answer is affirmative and, as shown in \cite{Gne(05)}, they coincide with the so--called $\sigma$--stable Poisson--Kingman models, which are obtained by means of a particular transformation of $\sigma$--stable completely random measures. The technical background needed for precisely defining such models goes beyond the scope of this review and we refer the interested reader to \cite{P03,Gne(05)}. For our purposes it is enough to note that the derivation of posterior quantities in this setting represents a challenging issue, which has not found a satisfactory solution to date.

So far we have provided various motivations, of theoretical and practical relevance, for the use of Gibbs--type priors and the sub--classes discussed in this section provide a glimpse of the nice and simple structure they
  share. Nonetheless, we still need to shed some light on another
  distributional aspect which is important for assessing their
  suitability for nonparametric inference, namely their support.
As mentioned in the Section~\ref{sec1}, a large topological support is a
desirable property for a nonparametric prior since the essence of
being nonparametric can be associated to the fact of
assigning prior positive probability to as many \comillas{candidate
  models} as possible. When considering the weak
topological support, which is the most natural in this framework, it can be shown (see \cite{Deb12}) that ``genuinely nonparametric'' Gibbs--type priors comply with this requirement and have full weak support: in other terms, any weak neighborhood of any distribution in
$\mathbf{P}_\X$ will have a priori positive probability. Here by \comillas{genuinely nonparametric} we mean Gibbs--type priors whose realizations are discrete distributions for which the number of support points is not bounded. This essentially boils down to considering Gibbs--type priors either with $\sigma\geq 0$ or with $\sigma<0$ and unbounded support of the prior $\pi$ on the number of components in \eqref{eq:hier}. Such priors can be shown to possess the full weak support property, i.e. their topological support coincides with the space of probability measures whose support is included in the support of the prior guess $P^*$. In particular, if the support of $P^*$ coincides with $\X$, the support of $Q$ is the whole space $\mathbf{P}_\X$.

\bigskip

\section{Hierarchical mixture models based on Gibbs--type priors}\label{sec:mixtures}

As outlined in the Introduction an important application of discrete
random probability measures and, then, of Gibbs--type priors occurs within hierarchical mixture models of the type \eqref{eq:hierar}: this corresponds to assuming exchangeable data $(Y_i)_{i\ge 1}$ from a random density defined by
\begin{equation}
{\tilde f}(y)=\int_{\mathbb{X}} f(y\mid x) {\tilde p}(\ddr x).
\label{eq:mixt_dens}
\end{equation}
In particular, when ${\tilde p}$ follows a discrete prior $Q$, a key ingredient for prior and posterior inferences is the corresponding EPPF. Indeed, given a set of observables $Y_1,\ldots,Y_n$ modeled according to the above random density, the clustering structure among the latent variables $X_1,\ldots,X_n$ drives both the posterior distribution on the number of components and the posterior density estimation. In particular,
\begin{equation*}
\pr(K_n=k\,\mid\, Y_1,\ldots,Y_n)
\propto \sum_{\mathsf{p}_n\in\mathcal{P}_{[n]}^k}
p_k^{(n)}(n_1,\ldots,n_k)\prod_{j=1}^{k}\int_{\mathbb{X}}
\prod_{i\in \mathcal{C}_j}f(y_i\mid x_j)P^*(\ddr x_j)
\label{eq:post_kn}
\end{equation*}
where $\mathcal{P}_{[n]}^k$ is the set of all partitions $\mathsf{p}_n$ of the $n$ latent variables into $k$ disjoint clusters and $\mathcal{C}_j$ identifies the indices of those latent variables $x_i$ that belong to the $j$--th cluster in the partition
$\mathsf{p}_n\in\mathcal{P}_{[n]}^k$. Therefore, the choice of $Q$ or, equivalently, of the corresponding EPPF, is crucial for
nonparametric Bayesian inferences in this framework and it can be
further appreciated through some numerical illustrations we are going
to provide later on in this section.

An appealing feature of Gibbs--type priors is their ability to control the prior mass allocated to different partitions through the reinforcement mechanism induced by the parameter  $\sigma$ and described in Section~\ref{sec:ramses}. This can be appreciated by looking at the induced (prior) distribution on the number $K_n$ of clusters. First note that the determination of the distribution of $K_n$ follows from a marginalization of \eqref{eq:1} and leads to
\begin{eqnarray}\label{pknGibbs}
\pr(K_n=k)=\frac{V_{n,k}}{\sigma^k}\,\mathscr{C}(n,k;\sigma)
\end{eqnarray}
with
\[
\mathscr{C}(n,k;\sigma)=\frac{1}{k!}\sum_{i=0}^k (-1)^i {k\choose i} (-i\sigma)_n
\]
denoting a generalized factorial coefficient. See \cite{Cha(05)} for
details on $\Ccr(n,k;\sigma)$. Substituting expressions
\eqref{vnkPD} and \eqref{vnkGG} in \eqref{pknGibbs} leads to the prior
distributions on the number of different components for the
PY and the NGG processes, respectively. Letting $\sigma\to 0$ in either of the resulting expressions one obtains the corresponding distribution for the  Dirichlet process case
    \begin{equation*}
    \pr(K_n=k)=\frac{\theta^k}{(\theta)_n}\, |s(n,k)|,
    \label{eq:Kn_Dirichlet}
    \end{equation*}
with $s(n,k)$ denoting the Stirling number of the first type. See \cite{Cha(05)}.

\begin{figure}[!htbp]
\begin{center}
\includegraphics[scale=0.6]{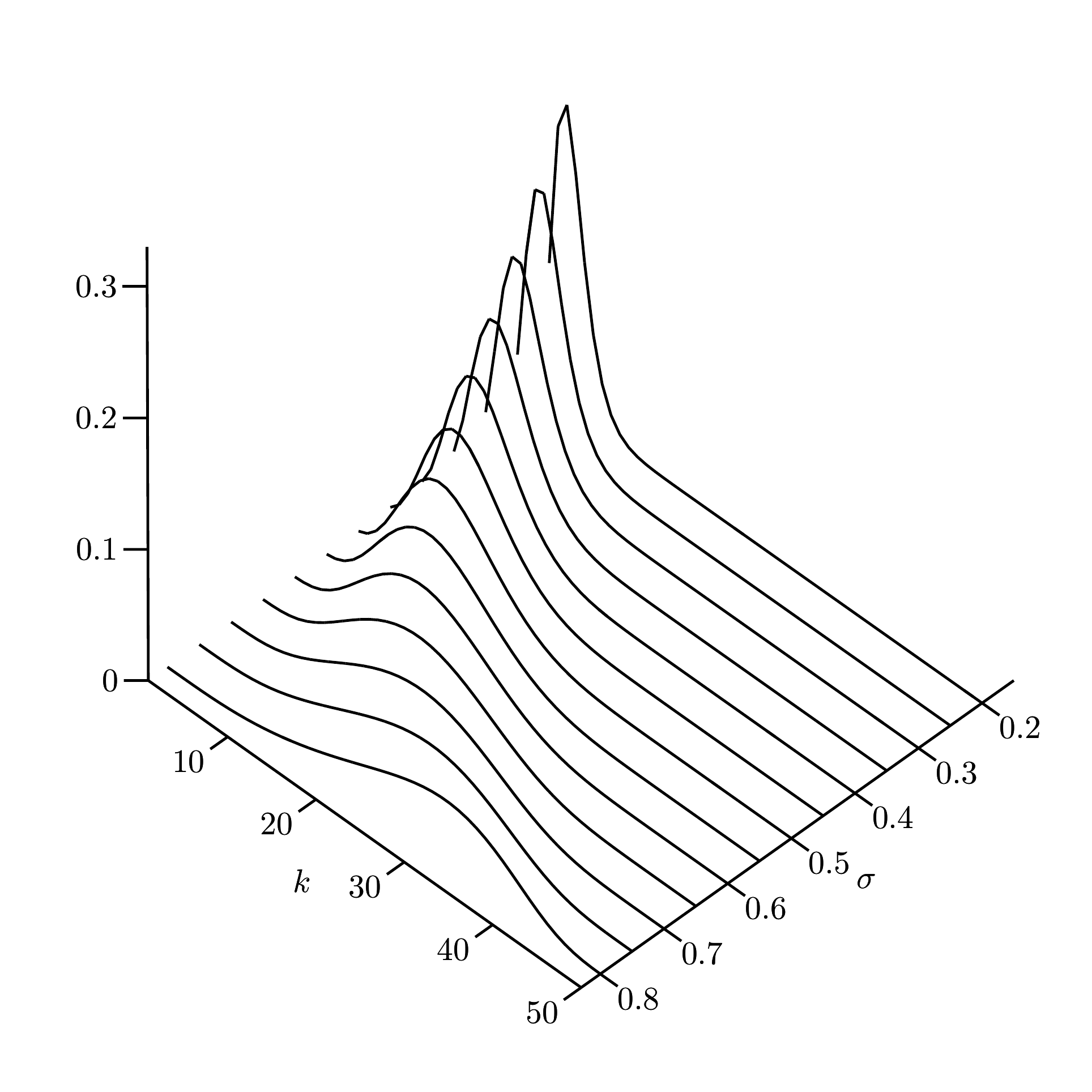}
\caption{Prior distributions on the number of groups corresponding to
  the NGG process with $n=50, \beta=1$ and $\sigma=0.2,0.3,\ldots,0.7$ and $\sigma=0.8$. }
\label{plt_pkgg_s}
\end{center}
\end{figure}

\begin{figure}[!htbp]
\begin{center}
\includegraphics[scale=0.8]{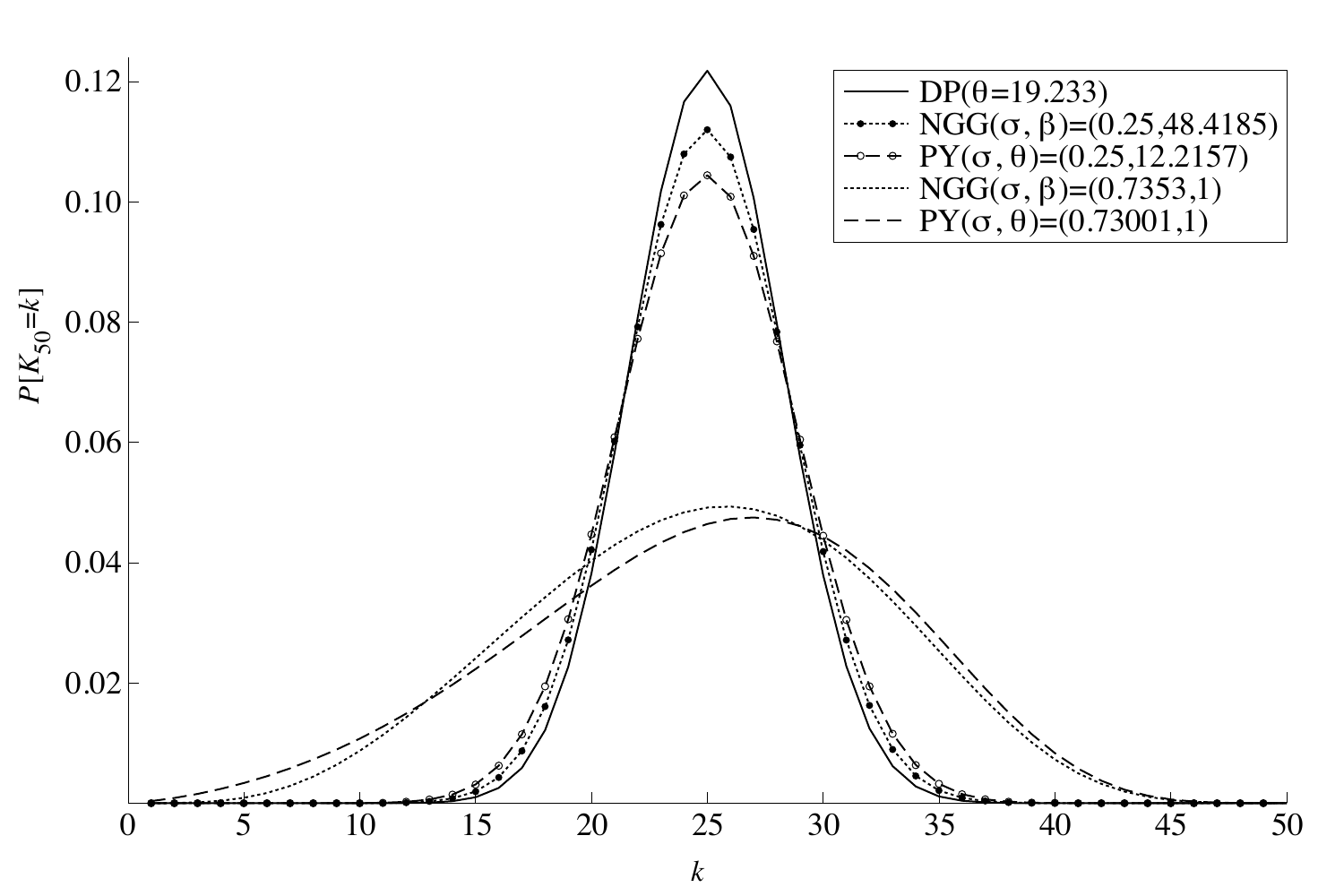}
\caption{Prior distributions on the number of clusters corresponding to  the Dirichlet (DP), the Pitman--Yor (PY) and the normalized generalized gamma (NGG) processes. The values of the parameters are set in such a way that $\E(K_{50})=25$.}
\label{pkn_com25}
\end{center}
\end{figure}

\begin{figure}[!htbp]
\begin{center}
\includegraphics[scale=0.8]{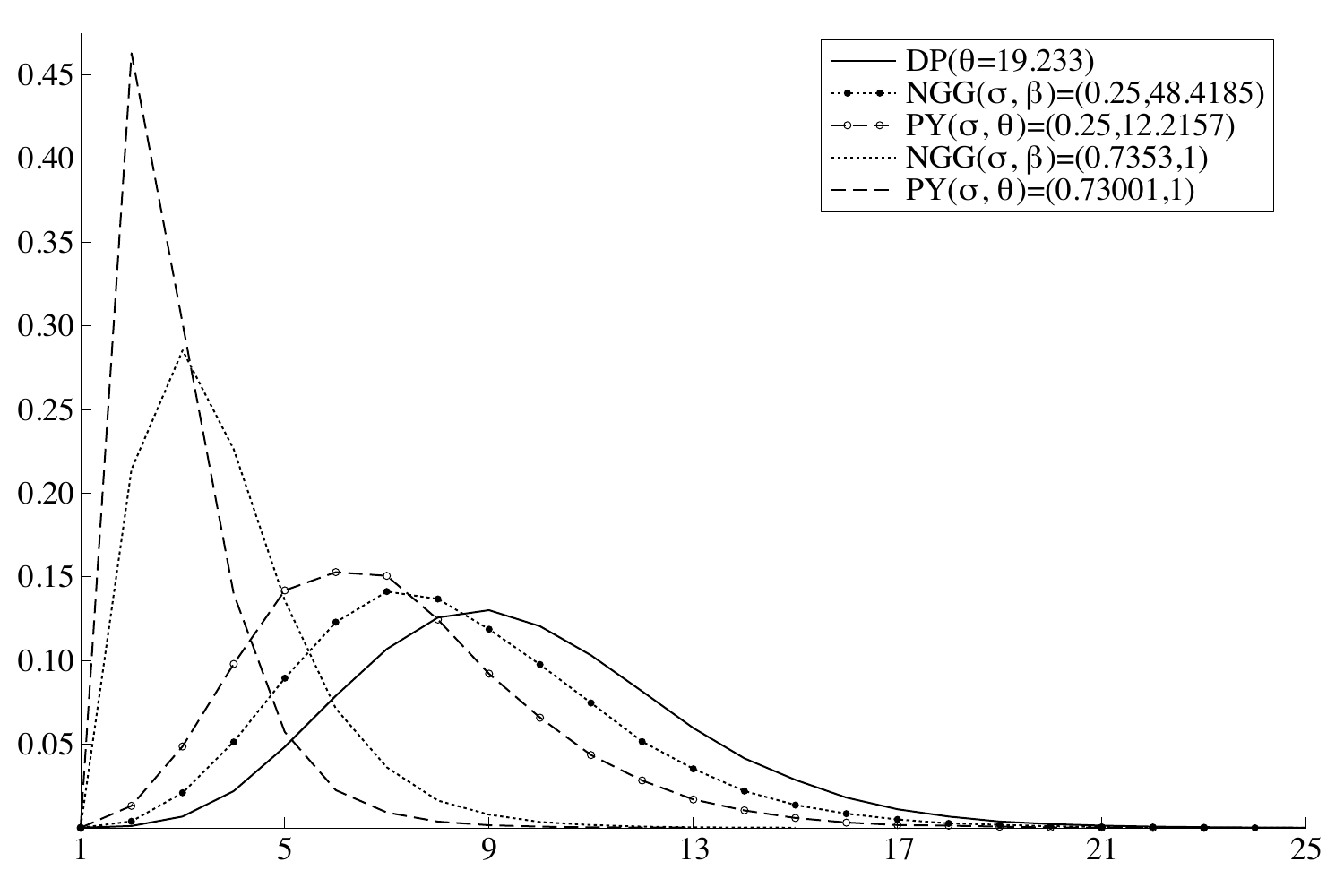}
\caption{Posterior distributions on the number of components corresponding to mixtures of the Dirichlet (DP), the Pitman--Yor (PY) and the normalized generalized gamma (NGG) processes with $n=50$ and parameters set so that $\E(K_{50})=25$.}

\label{postkn_com25}
\end{center}
\end{figure}

A graphical display of these distributions is best suited to highlight their differences. To this end, fix $n=50$ and consider the corresponding distributions of the number of components in the three above cases. For the Dirichlet process it is well--known that  the total mass parameter $\theta$ controls the location of the distribution of $K_{50}$: larger values of $\theta$ lead to a right-shift of the distribution implying an (a priori) larger number of components essentially affecting its dispersion. In both the PY process and NGG cases the role of controlling the location is played by $\theta$ and $\beta$, respectively. Hence, it is interesting to look at the additional parameter $\sigma$. Figure~\ref{plt_pkgg_s} concerns  the NGG process and displays the distribution of $K_{50}$ for a fixed value of $\beta$ and with $\sigma$ ranging between $0.2$ and $0.8$. Note that in Figures \ref{plt_pkgg_s}, \ref{pkn_com25} and \ref{postkn_com25} the probability masses are connected by straight lines only for visual simplification. From Figure \ref{plt_pkgg_s} it is evident that the addition of $\sigma$ allows to control the flatness, or the variability, of the distribution of $K_{50}$ thus yielding a higher degree of flexibility for the model. A similar behavior appears in the PY process.
Hence, replacing the Dirichlet process with a Gibbs--type prior characterized by a value of $\sigma$ in $(0,1)$ allows for a better control of the informativeness of the prior number of groups, since a larger $\sigma$ flattens the prior. To better visualize this fact, it is useful to consider a simple comparative example.
In addition to $n=50$, suppose that the prior expected number of clusters
is $25$. This implies that a reasonable criterion for eliciting the
parameters of a nonparametric prior is to fix them in a way such that
$\E(K_{50})=25$. We compare five different models: Dirichlet process with $\theta=19.233$, PY processes with
$(\sigma,\theta)=(0.25,12.2157)$ and $(\sigma,\theta)=(0.73001,1)$, and NGG processes with
$(\sigma,\beta)= (0.25,48.4185)$ and $(\sigma,\beta)=(0.7353,1)$, where all reported parameters are
chosen so that $E(K_{50})=25$. The corresponding distributions of
$K_{50}$ 
are displayed in Figure~\ref{pkn_com25}. Clearly, by increasing the value of $\sigma$ one obtains a less informative distribution on $K_{50}$: when moving from $\sigma=0$
to $\sigma\approx 0.73$ the distribution of $K_{50}$ becomes flatter, exhibiting a larger variability. The Dirichlet process, instead, implies a highly peaked distribution of $K_{50}$, which in terms of prior specification implies the need for a reliable prior information on the number of clusters, which is often unavailable. Furthermore, the PY and NGG processes have a similar behavior with the latter producing slightly lighter tails.

Let us now take a further step and compare the above five processes in
a toy example to have a closer look at the implication of such prior
specifications on posterior inferences on the clustering structure of
the data. To this end, assume that $n=50$ observations are drawn from a uniform mixture of two well-separated Gaussian distributions, $\mathsf{N}(1,0.2)$ and $\mathsf{N}(10,0.2)$. From a classification perspective these data clearly identify two groups.  We model them with the following nonparametric mixture model with standard specification
    \begin{align*}
    (Y_i\mid m_i, v_i)&\simind\mathsf{N}(
    m_i, v_i), \quad\qquad i=1,\ldots,n\\
    (m_i, v_i\mid {\tilde p})&\simiid {\tilde
      p}\qquad\qquad\qquad\quad i=1,\ldots,n\\
    {\tilde p}&\sim Q
    \end{align*}
with $Q$ corresponding to the five processes above and $P^*(\ddr m
,\ddr  v)=\mathsf{N}(m\mid \mu, \tau v^{-1})\mathsf{Ga}(v\mid 2,1)\,\ddr
m\, \ddr  v$, where $\mathsf{N}(\cdot\mid a, b)$
denotes the Gaussian density with mean $a$ and variance
$b>0$ and $\mathsf{Ga}(\cdot|c, d)$ is the density
corresponding to a Gamma distribution with mean $c/d$. A further
hierarchy is assumed for $\mu$ and $\tau$,
i.e. $\mu\sim\mathsf{N}(0,0.001)$ and
$\tau^{-1}\sim\mathsf{Ga}(1,100)$. In this setup the parameter
specification for the five processes (chosen so that
$E(K_{50})=25$) corresponds to a prior opinion on
$K_{50}$ remarkably far from the true number of components in the
mixture density that has generated the data. Given such a wrong prior
specification one then wonders whether the models possess enough
flexibility to shift a posteriori towards the correct number of
components, namely 2. The results are based on $100000$ iterations after $5000$ of burn in adopting a standard marginal MCMC algorithm with acceleration step. See \cite{EW95,mac} for further details on this algorithm. 

Figure~\ref{postkn_com25} depicts the posterior distribution on the
number of mixture components. The most important thing to note is that
a larger $\sigma$ leads to better posterior estimates. Both the
PY and NGG processes with $\sigma=0.73$, have been able to
shift most of the mass towards a very low number of components with
the PY process exhibiting a slightly better performance. See
also Table~\ref{Table:posttoy} for a display of the numerical values
of posterior probabilities associated to the possible values of
$K_{50}$.
This shows how a stronger reinforcement mechanism, implying a flatter
distribution of $K_{n}$, allows to recover more effectively the correct number of
components. In contrast, the Dirichlet process is
stuck around $10$ components, since the high peakedness of its prior
on $K_{n}$ prevents it from overruling completely the wrong prior
information. 

\begin{table}[htbp]
\begin{footnotesize}
\begin{center}
\begin{tabular}{|c|c|c|c|c|c|}
  \hline
$k$&DP($19.233$)& NGG($0.25$,$48.4185$) & PY($0.25$,$12.216$) & NGG($0.7353$,$1$) & PY($0.73001$,$1$)\\
  \hline
1	&0.0000	&0.0000	&0.0000	&0.0000	&0.0000	 \\
2	&0.0011	&0.0039	&0.0132	&0.2143	&{\bf 0.4630}	 \\
3	&0.0068	&0.0209	&0.0487	&{\bf 0.2854}	&0.3015	 \\
4	&0.0220	&0.0514	&0.0979	&0.2263	&0.1399	 \\
5	&0.0484	&0.0894	&0.1419	&0.1360	&0.0573	 \\
6	&0.0789	&0.1229	&{\bf 0.1528}	&0.0713	&0.0225	 \\
7	&0.1069	&{\bf 0.1412}	&0.1506	&0.0361	&0.0092	 \\
8	&0.1257	 	&0.1368	&0.1245	&0.0163	&0.0037	 \\
9	&{\bf 0.1301}	&0.1187	&0.0921	&0.0079	&0.0016	 \\
10	&0.1205	&0.0976	&0.0659	&0.0035	&0.0007	 \\
11	&0.1031	&0.0746	&0.0435	&0.0017	&0.0003	 \\
12	&0.0816	&0.0516	&0.0283	&0.0007	&0.0002	 \\
13	&0.0597	&0.0353	&0.0170	&0.0003	&0.0001	 \\
$\geq 14$&0.1151	&0.0556	&0.0237	&0.0004	&0.0001	     \\
  \hline
\end{tabular}
\end{center}
\caption{\rm Posterior distributions on the number of components
arising from mixtures of the Dirichlet process (DP), the
normalized generalized Gamma (NGG) process and the Pitman-Yor (PY)
process centered such that the prior expected value of the
number of components is $25$ with the sample size $n=50$.}\label{Table:posttoy}
\end{footnotesize}
\end{table}

\begin{figure}[!htbp]
\begin{center}
\includegraphics[scale=0.8]{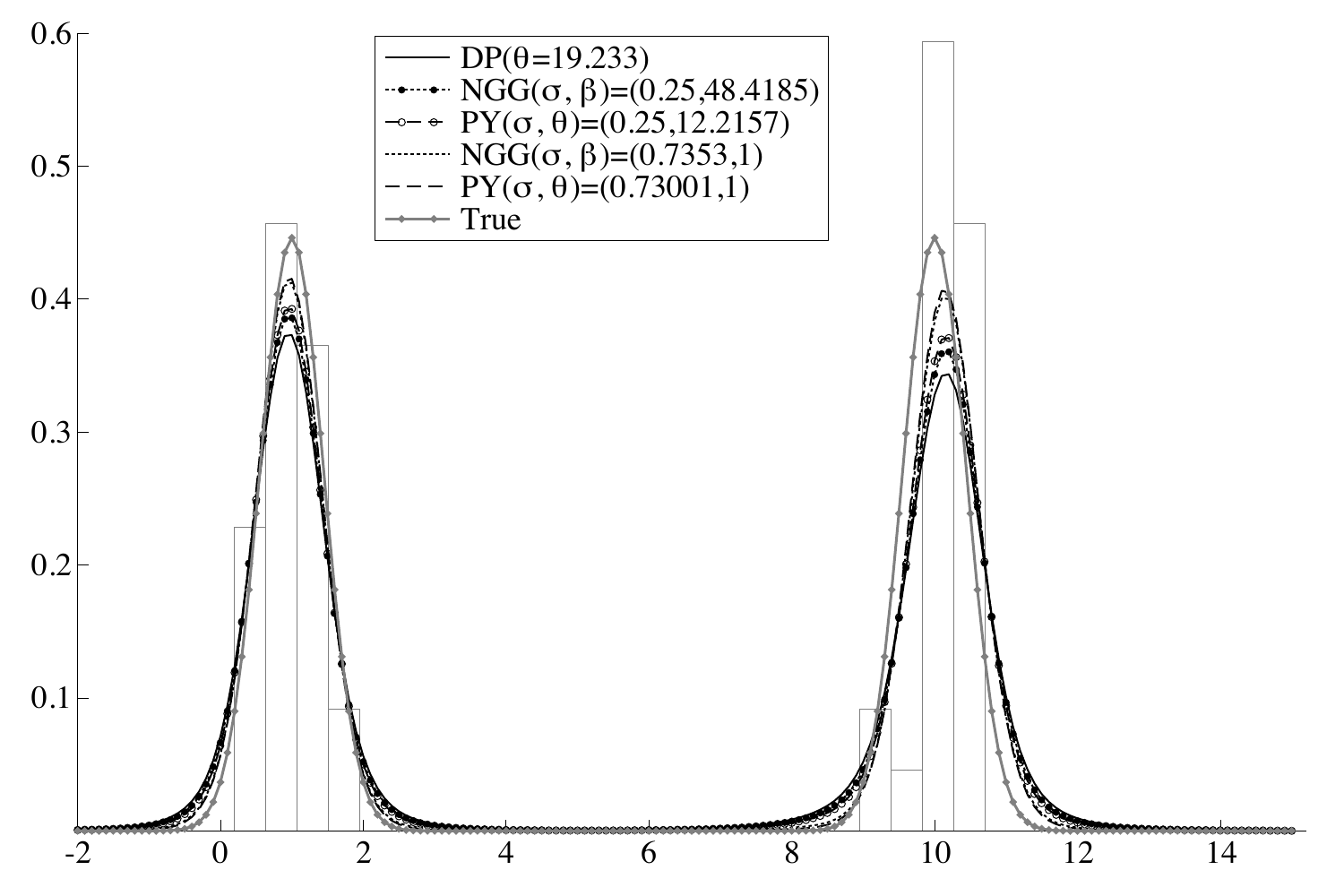}
\caption{Density estimates corresponding to the 5 mixture models that
  have been considered.}
\label{fig:density_post}
\end{center}
\end{figure}

Finally, it is important to point out that the above considerations
concerning the advantages of the additional parameter $\sigma$ hold beyond the present toy example since they represent structural
properties of the models, which are by now well--understood thanks to
several analytical results and computational analyses. See, e.g.,
\cite{LMP07}. As far as the estimates of the density $\tilde f$
  in \eqref{eq:mixt_dens} are concerned, these are displayed in
  in Figure~\ref{fig:density_post}. Even if the considerable
  heterogeneity in the posterior inferences on the number of
  components is not reflected by density estimates, one can still
  appreciate a slightly better performance of the NGG and PY processes
  with $\sigma=0.73001$ since they show a closer adherence
  to the depicted true density.

\bigskip

\section{Prediction in species sampling problems}\label{stefano}

As already mentioned in Section~\ref{sec1}, Gibbs--type priors are a powerful tool for addressing 
prediction and estimation 
in species sampling problems 
when observations are recorded from a population composed of individuals belonging to different types or species. This situation occurs in many applied research areas, including genetics, biology, ecology, economics and linguistics. Hence, in this section we will think of the observations $X_i$ in \eqref{eq:nonpmic} as species labels. The sample data $X_1,\ldots,X_n$ one can rely on for inferential purposes yield the
following pieces of information: the number $K_n$ of distinct species in the sample; the observed species labels $X_1^*,\ldots,X_{K_n}^*$; the frequencies $\bm{N}_n=(N_{1,n},\ldots,N_{K_n,n})$ associated to each of the observed species. Note that the last quantity can be alternatively reformulated in terms of the frequency counts $\bm{M}_n=(M_{1,n},\ldots,M_{n,n})$, where $M_{i,n}$ is the number of species that have appeared with frequency $i$ in the observed sample. It is obvious that these vectors must satisfy the following constraints:
    \[
    \sum_{i=1}^{K_n} N_{i,n}=n, \qquad
    \sum_{i=1}^n M_{i,n}=K_n,\qquad
    \sum_{i=1}^n iM_{i,n}=n.
    \]
In such problems species labels are typically not of interest, and the data can be efficiently summarized by either  $\bm{N}_n$ or $\bm{M}_n$,
namely the partition they form. Since the EPPF \eqref{eq:eppf_start} can
also be seen as the partition distribution induced by a sample, it
is natural to resort to the class of priors which have the most general yet tractable partition distribution. This naturally leads to work 
with Gibbs--type priors which are characterized by the product--form EPPF \eqref{eq:1}.

In this framework a novel Bayesian nonparametric methodology  for deriving estimators of quantities related to an
additional unobserved sample $X_{n+1},\ldots,X_{n+m}$ from $\tilde{p}$, conditional on $X_{1},\ldots,X_{n}$, has been proposed in \cite{Lij(07)} and \cite{lpw08}. An important applied problem is the estimation of the so-called overall \emph{species variety}, which can be measured by estimating the number $K_m^{(n)}=K_{n+m}-K_n$ of ``new'' distinct species that will be observed in the additional sample. A generalization has
been recently derived in \cite{Fav(12a)} and it corresponds to the estimator of the so-called \emph{rare species variety}
\begin{equation}
  \label{eq:rare}
  \widehat{M}_m^{(n)}(\tau)=
  \sum_{i=1}^\tau\widehat{M}_{i,m}^{(n)}=\sum_{i=1}^\tau\E(M_{i,n+m}\,|\,X_1,\ldots,X_n)
\end{equation}
namely the number of distinct species with frequency less than or
equal to a specific threshold of abundance $\tau$ that will be
detected in the additional sample of size $m$. Note
  that both the estimator of ${K}_m^{(n)}$, denoted by $\widehat{K}_m^{(n)}$, and
  $\widehat{M}_m^{(n)}(\tau)$ can be thought of as \emph{global measures} of
  overall and rare species variety respectively, since they are referred to
  the whole additional sample of size $m$. On the other hand, one
  may also need the corresponding \emph{local measures}, 
  which can be quantified in terms of the \emph{discovery probability} at step $(n+m+1)$ of the sampling process. 
  Bayesian estimators of the latter have been determined in \cite{Fav(12b)}. More specifically, if $\Delta_{i,n+m}$ is the set including species labels that appear with frequency $i\ge 0$ in the enlarged sample $X_1,\ldots,X_{n+m}$, one is interested in estimating
\begin{equation}
  \label{eq:discovery}
  U_{n+m,i}=\pr(X_{n+m+1}\in\Delta_{i,n+m}\,|\,X_1,\ldots,X_n).
\end{equation}
An estimator will be obtained by averaging over all possible realizations of the unobserved additional sample $X_{n+1}, \ldots, X_{n+m}$, conditional on the basic sample $X_1,\ldots,X_n$. Here $U_{n+m,0}$ stands for the probability of sampling a new species at step $(n+m+1)$, whereas $\sum_{i=0}^\tau U_{n+m,i}$ for the probability of sampling  either a species not yet observed or one with frequency less than $\tau$. Such local estimates are relevant, for example, in determining the size $m$ of the additional sample $X_{n+1},\ldots,X_{n+m}$: a possible criterion consists in fixing  $m$ equal to the maximum possible value for which the estimated discovery probability of new or rare species is above a certain threshold probability.

If the population is composed by a large number of unknown species (genes, agents, categories etc.) and the basic sample $X_1,\ldots,X_n$ displays only a small fraction of the species present in the population, Gibbs--type priors with $\sigma\in[0,1)$ are particularly suited. An effective and popular example is offered by the analysis of Expressed Sequence Tags \cite{Lij(07Bio)} or Serial Analysis of Gene Expression \cite{Gui(10)} data. Indeed, in these experiments either complementary DNA (cDNA) libraries or messenger RNA (mRNA) populations are considered and typical goals consist in identifying the genes they are composed of, the relative frequencies of such genes and also in comparing libraries/populations in terms of diversity. Due to time and cost constraints only a small portion of the whole library or population is typically sequenced and prediction is required to assess the overall characteristics. A similar experimental framework takes place in biological applications such as, for example, in the analysis of T-cell identification problems (see \cite{Sep(10)}). In this case one can characterize the immunological status of an organism by estimating the number of distinct clonotypes in a T-cell repertoire and the clonal size distribution, which is the frequency of clonotypes with a certain clonal size.
In contrast, if the population has a limited number of species, a common situation in Ecology, Gibbs--type models with $\sigma<0$ are more appropriate \cite{flmp}. In what follows, for brevity we will deal with the case $\sigma\in[0,1)$. This implies that, when specializing the results to the
PY process, one also has $\theta>-\sigma$. However, it is to be noted that most of the displayed findings carry over to the case of $\sigma<0$.


On this topic there exists a well--estalibished frequentist literature. The most relevant contributions typically draw inspiration from papers by I.J.~Good (\cite{Goo(53)}) and I.J.~Good and G.H.~Toulmin (\cite{Goo(56)}). See, e.g., \cite{Mao(02)} and \cite{Mao(04)}. For example, the popular Turing estimator for the discovery probability (displayed in \cite{Goo(53)} and credited to A.~Turing) is
    \begin{equation}
    \check{U}_{n,i}=(i+1)\frac{M_{i+1,n}}{n}.
    \label{eq:turing0}
    \end{equation}
For $i=0$ it provides an estimator for the probability that the $(n+1)$--th observation is new. 
Equivalently, $1-\check{U}_{n,0}$ provides an estimator of the sample coverage, namely the proportion of species observed in the sample, which is an important quantity in many applied frameworks. Moreover, estimates of $K_m^{(n)}$ and of the discovery probability $U_{n+m,0}$ for any $m\ge 1$ have been established in \cite{Goo(56)} and shall be henceforth termed Good--Toulmin estimators. They coincide with
    \begin{equation}
    \check{U}_{n+m,0}=n^{-1}\,\sum_{i=1}^\infty\,(-\lambda)^{i-1}\,i\:M_{i,n},\qquad
   \check{K}_m^{(n)}=\sum_{i=1}^\infty(-1)^{i-1}\lambda^i\,M_{i,n}
   \label{eq:good_toulmin}
    \end{equation}
where $\lambda=m/n$. Due to the alternating sign of the sums, when $\lambda$ is large they can yield inadmissible numerical values. This instability arises even for values of $m$ moderately large with respect to $n$, typically $m$ greater than $n$ is enough for it to appear. An illustration is provided in Section~\ref{sec:overall}. On the other hand, we are not aware of frequentist estimators of the discovery probabilities $U_{n+m,i}$ when both $m$ and $i$ are positive integers.

\subsection{Bayesian inference on overall species variety}\label{sec:overall}

Based on the EPPF~\eqref{eq:1}, an explicit expression for the distribution of the number of ``new'' distinct species observed in the additional sample, $K_{m}^{(n)}$, conditional on the information provided by $X_1,\ldots,X_n$, has been determined in \cite{Lij(07)} and is given by
    \begin{equation}\label{eq:estim_dist}
    \pr(K_{m}^{(n)}=j\,|\,X_1,\ldots,X_n)=
    \frac{V_{n+m,k+j}}{V_{n,k}}\frac{\mathscr{C}(m,j;\sigma,-n+k\sigma)}{\sigma^{j}}
    \end{equation}
where $X_1,\ldots,X_n$ is partitioned into $K_n=k$ clusters with respective frequencies $n_1,\ldots,n_k$ and $\Ccr(n,k;\sigma,\gamma)$ is the
non--central generalized factorial coefficient $$\mathscr{C}(m,j;\sigma,-n+k \sigma)=(j!)^{-1}\:\sum_{r=0}^j(-1)^{r}\,\binom{j}{r}\:(n-\sigma(r+k))_m.$$ See \cite{Cha(05)}. An important implication of \eqref{eq:estim_dist} is the sufficiency of $K_{n}$ for predicting the number of ``new'' distinct species. The expression in \eqref{eq:estim_dist} serves then as a basis for determining the Bayesian nonparametric estimator, with respect to a squared loss function, of the overall species variety as
    \begin{equation}\label{eq:est_gibbs}
    \widehat{K}_{m}^{(n)}=\text{E}(K_{m}^{(n)}\,|\, K_{n}=k,\: \bm{N}_{n}=\bm{n})
    \end{equation}
with $\bm{n}=(n_1,\ldots,n_k)$. This can be seen as a Bayesian counterpart of the Good--Toulmin estimator.

When $\tilde{p}$ is the PY process, with parameter
$(\sigma,\theta)$, \eqref{eq:estim_dist} becomes
\begin{equation}
\label{eq:estim_dist_twopar}
\pr(K_{m}^{(n)}=j\,|\,K_{n}=k,\bm{N}_{n}=\bm{n})= \frac{(\theta/\sigma+k)_{j}}{(\theta+n)_{m}}\mathscr{C}(m,j;\sigma,-n+k\sigma).
\end{equation}
As shown in \cite{Fav(09)}, the estimator for $K_{m}^{(n)}$ in \eqref{eq:est_gibbs} then reduces to
\begin{equation}\label{eq:est_twopar}
\widehat{K}_{m}^{(n)}=\left(k+\frac{\theta}{\sigma}\right)\left(\frac{(\theta+n+\sigma)_{m}}{(\theta+n)_{m}}-1\right).
\end{equation}
The main advantage of \eqref{eq:est_twopar}, and of other estimators
devised in \cite{Fav(09)} for measuring the overall species variety, is that they are explicit and can be exactly evaluated even when the size $m$ of the additional sample is large compared to the size of the basic sample $n$. This happens, for instance, in genomic applications where one has to deal with relevant portions of cDNA libraries consisting of millions of genes.

In several applied contexts it is useful to accompany point estimates
such as \eqref{eq:est_twopar} with the corresponding credible
intervals. These can be easily derived from the conditional distribution \eqref{eq:estim_dist_twopar}. However, if the sample sizes are very large the computation of the non-central generalized factorial coefficient may become cumbersome. To circumvent such a problem one could resort to asymptotic credible intervals. This motivates, also from a practical point of view, the asymptotic analysis of $K_{m}^{(n)}$, conditional on $K_{n}$, for a fixed $n$ and as $m\rightarrow\infty$, provided in \cite{Fav(09)}. Let $f_{\sigma}$ stand for the density function of a positive $\sigma$-stable random variable, and let $U_{q}$, for any $q\geq0$, be a positive random variable characterized by the following density function
\begin{displaymath}
f_{U_{q}}(u)=\frac{\Gamma(q\sigma+1)}{\sigma\Gamma(q+1)}u^{q-1-1/\sigma}f_{
\sigma}\left(u^{-1/\sigma}\right).
\end{displaymath}
Moreover, set $B_{a,b}$ as a beta random variable with parameters
$(a,b)$.
Conditional on the information provided by $X_1,\ldots,X_n$, one has
\begin{equation}\label{eq:limit_distinct}
\frac{K_{m}^{(n)}}{m^{\sigma}}\:\stackrel{\text{a.s.}}{\longrightarrow}\: Z_{n,k},
\end{equation}
as $m\rightarrow\infty$, where $Z_{n,k}\stackrel{\text{d}}{=}B_{k+\theta/\sigma,n/\sigma-k}U_{(\theta+n)/\sigma}$ with $B_{k+\theta/\sigma,n/\sigma-k}$ and $U_{(\theta+n)/\sigma}$ being independent. A similar asymptotic result has been obtained also
for the NGG process in \cite{Fav(12)}. Note that the
$\sigma$-diversity discussed in \eqref{eq:kappan_asympt} can be
recovered from \eqref{eq:limit_distinct} by setting $n=k=0$. Turning
back to the practical uses of \eqref{eq:limit_distinct} for the
determination of credible asymptotic intervals for $K_m^{(n)}$, it is
apparent that one still needs to derive the quantiles  of
$Z_{n,k}$. From an analytical point of view this is a challenging
task, which nonetheless can be avoided by 
resorting to straightforward computational algorithms that allow to sample from
the limit random variable $Z_{n,k}$, and thus to approximate the
quantiles. See \cite{Fav(09),Dev(09)} for details.

\subsection{Bayesian inference on rare species variety}

The problem of deriving estimators for the rare species variety has
been recently considered in \cite{Fav(12)} and \cite{Fav(12b)}. One of
such estimators is represented by the number of distinct species with
frequencies less than or equal to a specified threshold of abundance
$\tau$, for any $\tau\le n+m$, that are generated by the additional sample, as displayed in \eqref{eq:rare}. The determination of
$\widehat{M}_{i,m}^{(n)}$, under a square loss function,
is eased by resorting to the decomposition
\[
\widehat{M}_{i,m}^{(n)}=\widehat{N}_{i,m}^{(n)}+\widehat{O}_{i,m}^{(n)}
\]
where $\widehat{N}_{i,m}^{(n)}$ is the estimator of the number of
``new'' distinct species with frequency $i$ not detected in
$X_1,\ldots,X_n$ and $\widehat{O}_{i,m}^{(n)}$ is the estimator of the
number of ``old'' distinct species (i.e. included in $X_1,\ldots,X_n$) that appear with frequency $i$ in the enlarged sample. This implies that $\widehat{M}_m^{(n)}(\tau)$ in \eqref{eq:rare} arises as the sum of two well-defined quantities: (i) the estimator of the number of ``new'' distinct species with frequencies less than or equal to $\tau\leq m$ and generated by the additional sample, i.e. $\widehat{N}_m^{(n)}(\tau)=\sum_{i=1}^{\tau}\widehat{N}_{i,m}^{(n)}$; (ii) the estimator of the number of ``old'' distinct species with frequencies less than or equal to $\tau\leq n+m$ and generated by updating the frequencies of the partition induced by the basic sample with the additional sample,
i.e. $\widehat{O}_m ^{(n)}(\tau):=\sum_{i=1}^{\tau}\widehat{O}_{i,m}^{(n)}$. It is apparent that if $\tau=m$ one obtains $\widehat{N}_m ^{(n)}(\tau)=\widehat{K}_m^{(n)}$. In this respect, the concept of rare species variety can be interpreted as a generalization of the concept of overall species variety.

A result in \cite{Fav(12)} gives explicit expressions of the moments, of any order, of both the number of ``new'' species with frequency $i$ in $X_{n+1},\ldots,X_{n+m}$ and of the number of ``old'' species with frequency $i$ in the enlarged sample $X_1,\ldots,X_{n+m}$. From these one deduces $\widehat{N}_{i,m}^{(n)}$ and $\widehat{O}_{i,m}^{(n)}$ thus obtaining an estimator of rare species variety. It can be seen that
    \begin{equation}
    \label{eq:bayes_estimator_old}
    \widehat{O}_{i,m}^{(n)}=\sum_{t=1}^{i}{m\choose i-t}\, M_{t,n}\:(t-\sigma)_{i-t}\:
    \sum_{j=0}^{m}\frac{V_{n+m,k+j}}{V_{n,k}}\, \frac{\mathscr{C}\left(m-(i-t),j;\sigma,-n+t+(k-1)\sigma\right)}{\sigma^{j}}.
    \end{equation}
From \eqref{eq:bayes_estimator_old}, it is clear that  $(K_n,M_{1,n},\ldots,M_{\tau,n})$ is a sufficient statistic for predicting the number of ``old'' distinct species with frequency less than or equal to $\tau$. Moreover,
    \begin{equation}
    \label{eq:bayes_estimator_new}
    \widehat{N}_{i,m}^{(n)}= {m\choose i}(1-\sigma)_{i-1}
    \:\sum_{j=0}^{i}\frac{V_{n+m,k+j+1}}{V_{n,k}} \,\frac{\mathscr{C}(m-i,j;\sigma,-n+k\sigma)}{\sigma^{j}}
    \end{equation}
and $K_{n}$ is sufficient for predicting the number of ``new'' distinct species with frequency less than or equal to $\tau$. Finally, $\widehat{M}_{i,m}^{(n)}$ can be derived as the sum of the estimators in \eqref{eq:bayes_estimator_old} and \eqref{eq:bayes_estimator_new} and, then, $\widehat{M}_m^{(n)}(\tau)$ from \eqref{eq:rare}.

If we focus on the special case where the Gibbs--type prior is the PY process, then the expressions \eqref{eq:bayes_estimator_old} and \eqref{eq:bayes_estimator_new} considerably simplify and reduce to
    \begin{align*}
    \widehat{O}^{(n)}_{i,m}
    &=\sum_{t=1}^{i}{m\choose i-t}\, M_{t,n}\, (t-\sigma)_{i-t}\, \frac{(\theta+n-t+\sigma)_{m-(i-t) }}{(\theta+n)_{m }}\\[4pt]
    \widehat{N}^{(n)}_{i,m}
    &={m\choose i}(1-\sigma)_{i-1 }(\theta+k\sigma)\, \frac{(\theta+n+\sigma)_{m-i }}{(\theta+n)_{m}}.
    \end{align*}
It is worth noting that the determination of estimators of the rare species variety poses a major technical hurdle that does not occur when estimating the overall species variety. Indeed, one has to consider all possible modifications, induced by the observations in the additional sample, on the frequencies of the species detected in the basic sample.

In the special PY process case, one can establish the asymptotic behavior of rare species variety as $m\to\infty$. This is somehow in the spirit of \eqref{eq:limit_distinct} in the context of overall species variety. In this case, as shown in \cite{Fav(12)}, one has for any $i\ge 1$
    \begin{displaymath}
      \frac{M_{i,n+m}\,|\,X_1,\ldots,X_n}
      {m^{\sigma}}\stackrel{d}{\longrightarrow}
    \frac{\sigma(1-\sigma)_{i-1}}{i!}\: Z_{n,k},
    \end{displaymath}
as $m\rightarrow\infty$, where $Z_{n,k}$ is the limit random variable introduced in \eqref{eq:limit_distinct} and $\stackrel{d}{\longrightarrow}$ stands for convergence in distribution. This implies that $K_{n}$ is asymptotically sufficient for predicting the number of distinct species with frequency $i$ that are generated after observing the additional sample, conditional on the information provided by the random partition of the basic sample.

Rare species variety can be further assessed locally in terms of
discovery probabilities $U_{n+m,i}$ as defined in
\eqref{eq:discovery}. This leads to the proposal of Bayesian nonparametric counterparts to the Turing and the Good--Toulmin estimators that are recalled in \eqref{eq:turing0} and in \eqref{eq:good_toulmin}, respectively. If one assumes a square loss function, then an estimator of $U_{n,i}$ is
    \begin{equation}
    \label{eq:bnp_turing}
    \widehat{U}_{n,i}=\frac{V_{n+1,k}}{V_{n,k}}\:(i-\sigma)\:M_{i,n}
    \end{equation}
for any $i\le n$, while the discovery probability of a new species, i.e. $i=0$, can be easily deduced from the predictive distribution \eqref{eq:predict_gibbs} and is given by $\widehat{U}_{n,0}=V_{n+1,k+1}/V_{n,k}$. Note that, unlike the Turing estimator, $\widehat{U}_{n,i}$ depends on $M_{i,n}$ which seems to be more coherent with what intuition would suggest. If we now let $m\ge 1$ and $j\le n+m$, an estimator of the discovery probability turns out to be\footnote{The estimators in \eqref{eq:estimated_abundance_coverage_probability} and \eqref{eq:estimated_abundance_coverage_probability_pd} slightly differ from those in \cite{Fav(12b)}, since the latter contain a minor inaccuracy that we have corrected here.}
\begin{align}
\label{eq:estimated_abundance_coverage_probability}
\widehat{U}_{n+m,i}&=\sum_{l=1}^{i}M_{l,n}(l-\sigma)_{i+1-l}{m\choose
  i-l}Q_{m,i}^{(n,k)}(l,0,l-\sigma)\\[4pt]
&\qquad\qquad + \sigma(1-\sigma)_{i}{m\choose i}Q_{m,i}^{(n,k)}(1,1,0)
\notag
\end{align}
where
\begin{displaymath}
Q_{m,i}^{(n,k)}(\alpha,\beta,\gamma)=\sum_{r=\beta}^{m-i+\alpha} \frac{V_{n+m+1,k+r}}{V_{n,k}}\frac{\mathscr{C}(m-i+\alpha-\beta,r-\beta;\sigma,-n+k\sigma+\gamma)}{\sigma^{r}}.
\end{displaymath}
When $i=0$ and $m\ge 1$, this yields 
the following Bayesian analog of the Good--Toulmin estimator for the probability of discovering a new species 
\begin{equation}
  \label{eq:bnp_goodtoulmin}
  \widehat{U}_{n+m,0}=\sum_{i=0}^m \frac{V_{n+m+1,k+i+1}}{V_{n,k}}\:
  \frac{\Ccr(m,i;\sigma,-n+k\sigma)}{\sigma^i}
\end{equation}
whereas there is no frequentist counterpart to \eqref{eq:estimated_abundance_coverage_probability} when both $m$ and $k$ are positive integers. From these closed form expressions one can deduce a further measure of rare species variety as $\widehat{U}_{n+m}(\tau)=\sum_{i=0}^\tau \widehat{U}_{n+m,i}.$

If one adopts a specification of the $V_{n,k}$'s yielding a PY process, nice and simple forms of the estimators of the discovery probabilities and of rare species variety are obtained. For example, the analog of the Turing estimator \eqref{eq:turing0}
reduces to
\begin{equation*}
  \label{eq:bnp_turing_py}
  \widehat{U}_{n,i}=\frac{i-\sigma}{\theta+n}\:M_{i,n}
\end{equation*}
and the Bayesian counterpart \eqref{eq:bnp_goodtoulmin} of the Good--Toulmin estimator coincides with
\begin{equation}
  \label{eq:bnp_goodtoulm_py}
  \widehat{U}_{n+m,0}=\frac{\theta+k\sigma}{\theta+n}\:
  \frac{(\theta+n+\sigma)_m}{(\theta+n+1)_m}.
\end{equation}
Finally, the probability that the $(n+m+1)$--th observation coincides with a species detected $j$ times in the enlarged sample
$X_1,\ldots,X_{n+m}$ is
    \begin{equation}\label{eq:estimated_abundance_coverage_probability_pd}
      \widehat{U}_{n+m,i}
      =\sum_{l=1}^{i}M_{l,n}(l-\sigma)_{i+1-l}{m\choose
        i-l}\frac{(\theta+n-l+\sigma)_{m-i+l}}{(\theta+n)_{m+1}}
      + (1-\sigma)_{i}{m\choose i} \frac{(\theta+k\sigma)(\theta+n+\sigma)_{m-i}}{(\theta+n)_{m+1}}.
    \end{equation}
To briefly illustrate the behavior of the Bayesian nonparametric estimator based on the PY process and compare it with the Good--Toulmin estimator let us consider genomic data, which consists of Expressed Sequence Tags (EST) obtained from \emph{Naegleria gruberi} cDNA libraries. \emph{Naegleria gruberi} is a widespread free--living soil and freshwater \emph{amoeboflagellate} widely studied in the biological literature. The two considered datasets are sequenced from two cDNA libraries prepared from cells grown under different culture conditions, aerobic and anaerobic, and have been previously analyzed in \cite{Susko2004,Lij(07Bio)}. The sequenced data, which will constitute the basic samples, are reported in Table~\ref{table:counts}.

\begin{table}
\begin{footnotesize}
\label{table:counts}
\begin{center}
  \begin{tabular}{|l|c|c|c|c|c|c|c|c|c|c|c|}
    \hline
    {\it Library }
    &{\it 1}&{\it 2}&{\it 3}&{\it 4}&{\it 5}&{\it 6}&
    {\it 7}&{\it 8}&{\it 9}&{\it 10}&{\it 11}\\
    \hline
    {Naegleria Aerobic}          &346&57&19&12&9 &5&4&2&4&5 &4 \\ 
    {Naegleria Anaerobic}        &491&72&30&9 &13&5&3&1&2&0 &1 \\
    \hline
    {\it Library}
    &{\it 12}&
    {\it 13}&{\it 14}&{\it 15}&{\it 16}&{\it 17}&{\it 18}&
    {\it 27}&{\it 55}& $j$ & $n$ \\ \hline
    {Naegleria Aerobic}          &1 &0 &0 &0 &1 &1 &1 &1 &1 &473&959\\
    {Naegleria Anaerobic}        &0 &1 &3 &0 &0 &0 &0 &0 &0 &631&969\\
    \hline
  \end{tabular}
\end{center}
\caption{\rm ESTs from two {\it Naegleria gruberi} libraries. Reported data include: frequency counts $M_{i}$, for different values of $i$, total number of distinct genes $j$ and sample size $n$. Source: Susko and Roger (2004).}
\end{footnotesize}
\end{table}

\begin{figure}[!htbp]
\begin{center}
\includegraphics[scale=0.8]{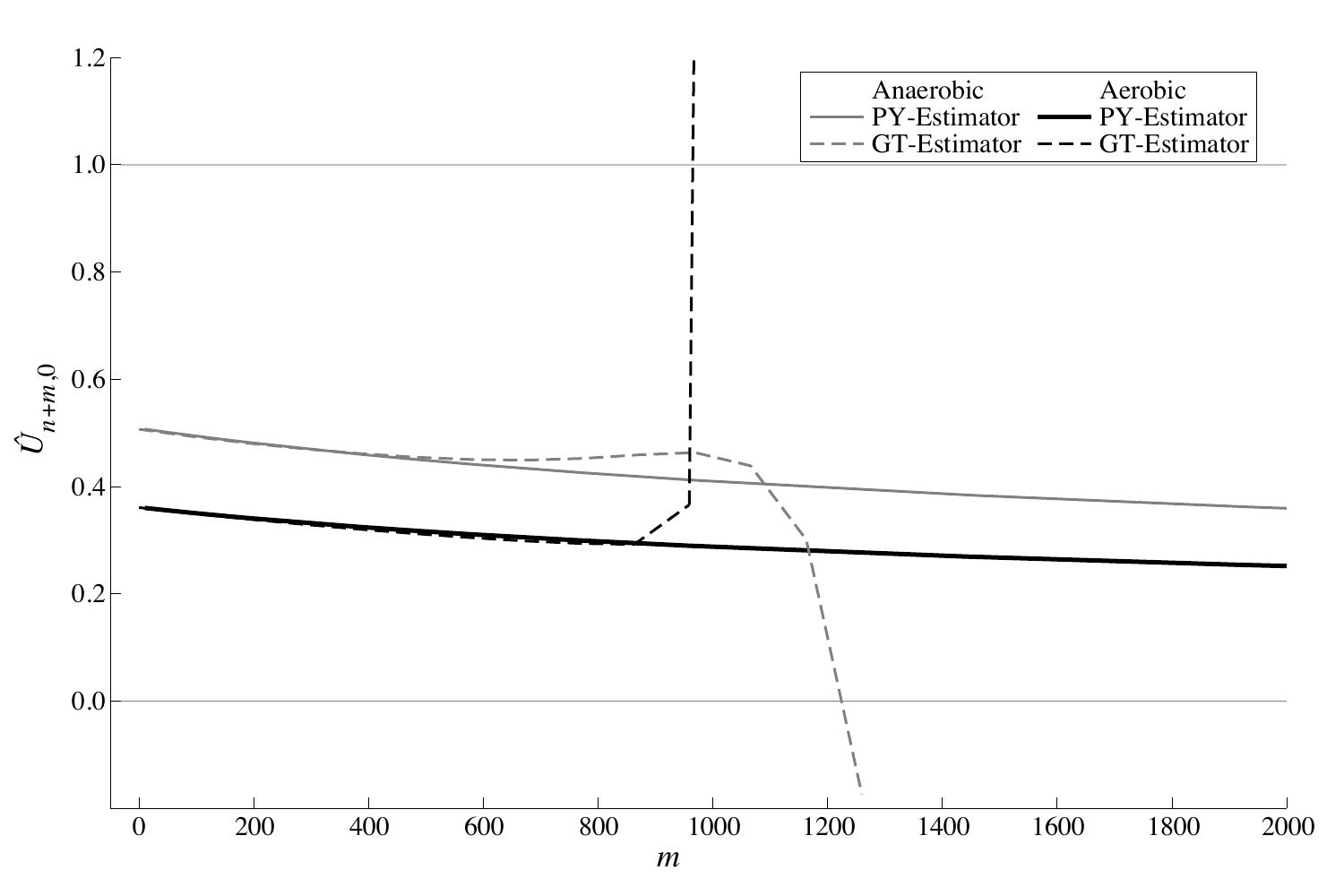}
\caption{EST data from \emph{Naegleria gruberi} aerobic and anaerobic cDNA libraries with basic sample $n \cong 950$: Good--Toulmin (GT) and Pitman--Yor process (PY) estimators of the probability of discovering a new gene at the $(n+m+1)$--th sampling step for $m=1, \ldots, 2000$.}\label{fig:GT}
\end{center}
\end{figure}

If one is interested in the probability of discovering
a new gene at the $(n+m+1)$--th step of the sequencing process, one has two
options: the Good--Toulmin estimator $\check{U}_{n+m,0}$ reported in
\eqref{eq:good_toulmin} or the estimator $\hat{U}_{n+m,0}$ in \eqref{eq:bnp_goodtoulm_py}  which is based on the PY process. To
complete the specification of the latter let us mention that the
parameters $(\sigma,\theta)$ are fixed according to an empirical Bayes
specification, which yields $(0.66,155.5)$. The results are displayed in Figure~\ref{fig:GT}. It is clear that the Good--Toulmin estimator exhibits an erratic behavior for values of the additional sample 
relatively larger than that of the basic sample $n$. 
This phenomenon 
is avoided by the Bayesian nonparametric estimator since it relies on a well--defined probabilistic model in which all quantities are modeled jointly and coherently. For sizes of $m$ for which the Good--Toulmin estimator works well, the estimators essentially coincide. Note that, in terms of the specific application, the anearobic library exhibits the clearly higher genetic diversity. Furthermore, as already mentioned, one can use such estimates to fix the size of the additional sample $m$ as the maximum integer for which the discovery probability lies above the desired threshold, which is typically determined also on the basis of cost considerations.


\bigskip

\section{Frequentist asymptotics}\label{sec:consistency}

During the last two decades frequentist consistency has gained a major role in Bayesian nonparametrics and is generally accepted as a key validation criterion for the use of a nonparametric prior in applied problems.  See \cite{Gho10} for a recent review on the subject. The idea that underlies the study of consistency consists in assuming that the data $(X_n)_{n\geq 1}$ are iid from some \comillas{true} $P_0\in\mathbf{P}_X$ and in verifying whether the posterior distribution $Q(\,\cdot\, | X_1,\ldots,X_n)$ accumulates in any suitably defined neighborhood of $P_0$. Therefore, while the posterior is derived based on an assumption of exchangeability of the data as described in \eqref{eq:nonpmic}, the frequentist asymptotic evaluation postulates plain independence of the data generating process. This explains why such an approach has also been termed \comillas{what if} approach by P. Diaconis. See \cite{Dia86}. Here we shall discuss consistency for Gibbs--type priors. In this respect, note that frequentist asymptotics of Bayesian procedures is different from the kind of asymptotics discussed in Section~\ref{stefano} which preserves a Bayesian flavor since it aims at achieving a large sample approximation of the posterior without modifying the dependence assumption among the data.

Let us start by fixing some notation and introducing some useful concepts. First the data $X_i$'s are assumed to be iid from some \comillas{true} $P_0$ or, in other terms, the distribution of the sequence of observations $(X_n)_{n\ge 1}$ is the infinite product measure $P_0^\infty=P_0\times P_0\times\cdots$. If $A_\ep$ denotes a neighborhood of $P_0$ of radius $\ep$, the posterior is said to be consistent at $P_0$ if
  $Q(A_\epsilon|X_1,\ldots,X_n)\to 1$
almost surely with respect to $P_0^\infty$, as $n\to \infty$ and for
any $\epsilon>0$. In the case of Gibbs--type priors, the natural choice for $A_\ep$ is represented by
weak neighborhoods. Clearly, consistency can be achieved
only at $P_0$ whose weak neighborhoods have a priori positive
probability. In this respect, the full support property of Gibbs--type
priors recalled in Section~\ref{sec:sub} is important since it ensures that
consistency can potentially be achieved at any given
$P_0$. Furthermore, note that the full support property represents a
desirable property not only when studying consistency in the case
where Gibbs--type priors are used to model directly the data, but also
in the context of mixture models as in \eqref{eq:hierar}. Indeed,
together with some other features of Gibbs--type priors, it allows to
extend known consistency results for Dirichlet process mixture models in a straightforward way and the condition for it to hold will be essentially the same. See \cite{Gho99,Lij05}.

As explained in some detail below, recent results
  suggest that Gibbs--type priors are always consistent
  with respect to (w.r.t.) discrete $P_0$'s. On the other hand, when they are used to
  model data coming from diffuse distributions, inconsistency may
  arise. 
  Possible inconsistency at diffuse $P_0$ should not,
  however, be interpreted as a serious issue: what really matters
  is the data generating mechanism the nonparametric prior is designed for so that
  consistency must hold w.r.t. choices of $P_0$ that are compatible
  with such a mechanism. Since Gibbs--type priors
  are discrete random probability measures, one should be primarily
  interested in investigating consistency w.r.t. discrete
  $P_0$'s. Indeed, Gibbs--type priors, and discrete nonparametric
  priors in general, are designed to model discrete
  distributions and they should under no circumstance be used to model
 data coming from
 diffuse distributions. In the latter case they should be exploited
 within hierarchical mixtures.

\subsection{General results}

The strategy for showing consistency consists in first identifying the
weak limit of the posterior, say $P'$ in $\mathbf{P}_\X$, which will
be some function of $P_0$, and then checking whether $P'=P_0$ so that
consistency is achieved.
The candidate weak limit $P'$ is identified by investigating the asymptotic behavior of the predictive distribution \eqref{eq:predict_gibbs} (i.e. the posterior expected value), which in explicit cases allows to guess $P'$ quite easily. Then one has to show that the posterior variance of $\p$ in \eqref{eq:discrete} converges to $0$, a.s.-$P_0^\infty$, which suffices to establish that the posterior concentrates in a weak-neighborhood of the predictive distribution. See \cite{Jam08,Deb12} for details.
Now, let $X_1,\ldots,X_n$ denote a sample with $\kappa_n$ distinct
values with corresponding frequencies
$n_1,\ldots,n_{\kappa_n}$. Even if $\kappa_n$ denotes
  the same quantity identified as $K_n$ in previous sections, we shall
  use a different symbol to emphasize the fact that here the asymptotic behavior of the
number of observed distinct species $\kappa_n$ is dictated by $P_0$
from which the iid sequence is sampled and not by a
  Gibbs--type prior directing an exchangeable sequence according to \eqref{eq:nonpmic}. Different choices of $P_0$ clearly yield different (almost sure) limiting behaviors for $\kappa_n$. On the one hand, if  $P_0$ is discrete with $N$ point masses, for any  $N\in\N\cup\{\infty\}$, then $P_0^\infty(\lim_n \kappa_n=N)=1$ and $P_0^\infty(\lim_n n^{-1}\kappa_n=0)=1$ even if $N=\infty$. On the other hand, if $P_0$ is diffuse, $P_0^\infty(\kappa_n=n)=1$ for any $n\ge 1$. Henceforth we shall focus on these two cases and adopt the shorter notation $\kappa_n\ll_{a.s.}n$ and $\kappa_n\sim_{a.s.}n$, which stand for $\kappa_n/n \to 0$ and $\kappa_n/n \to 1$ a.s.-$P_0^\infty$, respectively. 
It turns out that a key quantity for studying the asymptotics of the predictive distribution is given by the probability \eqref{eq:pred_new} of discovering a new observation at the $(n+1)$--th sampling step, which is given by $V_{n+1,\kappa_n+1}/V_{n,\kappa_n}$ in the case of Gibbs--type priors. Considering a Gibbs--type prior with base measure $P^*$ having support $\X$ and assuming that
\begin{equation}
  \frac{V_{n+1,\kappa_n+1}}{V_{n,\kappa_n}}\to \alpha\quad
  \mbox{a.s.-}P_0^{\infty}
  \tag{H}  \label{eq:A}
\end{equation}
as $n\to\infty$ for some constant $\alpha\in[0,1]$, in \cite{Deb12} it is shown that
\begin{equation*}
  \label{eq:weak}
  Q(A_\epsilon'|X_1,\ldots,X_n)
  \:\to\: 1 \qquad\qquad\mbox{a.s.-}P_0^\infty
    \end{equation*}
as $n\to\infty$ and for any $\epsilon>0$  where $A_\epsilon'$ is a weak neighborhood of $P'$. Moreover, one has
\begin{equation} \label{eq:linear_comb}
  P'=\alpha P^*(\cdot)+(1-\alpha)P_0(\cdot).
\end{equation}
Some comments regarding the above convergence result are in order. As
for the condition \eqref{eq:A}, it is worth noting that it holds true for all Gibbs--type priors for which an explicit expression of the $V_{n,\kappa_n}$'s is known, regardless as to whether $P_0$ is discrete or diffuse. It therefore represents only a mild regularity condition. Moreover, the posterior distribution converges to a point mass at \eqref{eq:linear_comb}, a linear combination of the prior guess $P^*$ and the ``true'' distribution $P_0$. Hence, weak consistency is guaranteed if $\alpha=0$ (and in the trivial case $P^*=P_0$ to be excluded henceforth) and it is sufficient to check whether the probability of discovering a new value converges to $0$, a.s.-$P_0^\infty$. Also, one can assess the departure from consistency by looking at the size of $\alpha$: the larger $\alpha$, the heavier the limiting mass assigned to the prior guess $P^*$. One can even think of a case of ``total inconsistency'', i.e. $\alpha=1$, the worst case scenario where the posterior tends to concentrate around the prior guess $P^*$ and no learning at all takes place.

To better visualize the above convergence result it is useful to look
at special cases of the PY process with $\sigma\in[0,1)$ and
$\theta>-\sigma$, for which such convergence had already been
established in \cite{Jam08}. From the form of their
predictive distributions \eqref{eq:pred_PY}, one can immediately
conjecture the following result: when $P_0$ is discrete ($\kappa_n\ll_{a.s.}n$) we have $\alpha=0$, implying consistency; when $P_0$ is diffuse ($\kappa_n\sim_{a.s.}n$), we have $\alpha=\sigma$, hence inconsistency, unless $\sigma=0$, which corresponds to the Dirichlet case. See also \cite{Jan10}. An analogous result has been established  for the NGG process together with some results concerning the case of Gibbs--type priors with $\sigma>0$ in \cite{Jan10}.

Focusing now on Gibbs--type priors with $\sigma<0$
  allows to highlight the occurrence of interesting phenomena. Recall from Section
  \ref{sec:sub} that these priors coincide with mixtures of
  PY processes with parameters $\{(\sigma,\, m|\sigma|):\:
  m=1,2,\ldots,\}$ and they can be represented in hierarchical form as
  \eqref{eq:hier}. It turns out that, according to the nature of the
  \comillas{true} distribution $P_0$, a sufficient condition can be
  stated in terms of the tail behavior of the mixing distribution $\pi$ in \eqref{eq:hier}. More precisely, for Gibbs--type priors with parameter $\sigma<0$ and prior guess $P^*$ whose support coincides with $\X$, in \cite{Deb12} consistency is shown to hold
\begin{itemize}
\item[{\rm (i)}] at any discrete $P_0$ if for sufficiently large $m$
\begin{equation}
  \frac{\pi(m+1)}{\pi(m)}\leq 1;
  \tag{T1}  \label{eq:B}
\end{equation}
\item[\rm{(ii)}] at any diffuse $P_0$ if for sufficiently large $m$ and for some $M<\infty$
\begin{equation}
  {\pi(m+1)\over \pi(m)}\leq {M\over m}.
  \tag{T2}  \label{eq:antonio2}
\end{equation}
\end{itemize}
Condition \eqref{eq:B} is an extremely mild assumption on the regularity of the tail of the mixing $\pi$: it requires $x\mapsto \pi(x)$ to be ultimately decreasing, a condition met by the commonly used probability measures on $\N$. Hence one can conclude that Gibbs--type priors with parameter $\sigma<0$ are essentially always consistent when $P_0$ is discrete. On the other hand, condition \eqref{eq:antonio2} requires the tail of $\pi$ to be sufficiently light, so when $P_0$ is diffuse one needs to closely investigate the tail behavior of $\pi$.

\subsection{Illustrations}\label{sec:cons_ex}

In light of the results stated above one is naturally led to wonder what happens when \eqref{eq:antonio2} is not satisfied. To this end we consider three different Gibbs--type priors presented in Section \ref{sec:sub} with $\sigma=-1$: each prior is characterized by a specific choice of the mixing distribution $\pi$. We focus on the case of diffuse $P_0$, which leads to some interesting conclusions. In the case of discrete $P_0$ it is straightforward to show that \eqref{eq:B} holds, hence ensuring consistency.

The first prior we consider, introduced in \cite{Gne10}, is characterized by the heavy-tailed mixing distribution \eqref{eq:mixing_gnedin}, which does not admit a finite expected value. Since $\pi(m+1)/\pi(m)=(m-\gamma)/(m+1)$ cannot be eventually bounded by $M/m$ for some constant $M$, condition \eqref{eq:antonio2} does not hold true. Given the $V_{n,\kappa_n}$'s admit the simple closed form expression \eqref{eq:gnedin_model}, the weights of the prediction rule simplify to
  $${V_{n+1,\kappa_n+1}\over V_{n,\kappa_n}}
  ={\kappa_n(\kappa_n-\gamma)\over n(\gamma+n)}.$$
It is easy to see that, if $P_0$ is diffuse, implying $\kappa_n\sim_{a.s.}n$, condition \eqref{eq:A} holds true with $\alpha=1$ and the weak limit coincides with the prior guess $P^*$, whatever the ``true'' distribution of the data $P_0$. This means we are in the case of \comillas{total} inconsistency.

The second example has a Poisson mixing distribution \eqref{eq:mixing_poisson} on the positive integers. Such a $\pi$ has light tails and condition~\eqref{eq:antonio2} is satisfied since $\pi(m+1)/\pi(m)=\lambda/(m+1)$. Therefore, by \eqref{eq:antonio2}, the posterior is consistent when $P_0$ is diffuse. 

The last sub-family of Gibbs--type priors with $\sigma=-1$ is identified by a geometric mixing distribution \eqref{eq:mixing_geometric}. Note that $\pi(m+1)/\pi(m)=\eta$ so that condition \eqref{eq:antonio2} does not hold true. It turns out that, with $P_0$ diffuse and $\kappa_n\sim_{a.s.}n$, one obtains
\begin{equation}\label{eq:asym_geom}
  {V_{n+1,\kappa_n+1}\over V_{n,\kappa_n}}
  \to\alpha={2-\eta-2\sqrt{1-\eta}\over \eta}\in[0,1].
\end{equation}
See \cite{Deb12} for details. The limit $\alpha$ in \eqref{eq:asym_geom} can be any point in $[0,1]$ according to the value of $\eta$ and therefore we can obtain the whole spectrum of weak limits \eqref{eq:linear_comb}
ranging from consistency ($\alpha=0$) to \comillas{total} inconsistency ($\alpha=1$). In particular, $\alpha$ is increasing in $\eta$, so the larger $\eta$, the heavier the limiting mass assigned to the prior guess. Small values of $\eta$ identify a situation similar to the second example since they yield a light-tailed $\pi$. Conversely, large values of $\eta$ are more in line with the first example giving rise to heavy-tailed $\pi$. Finally, it is worth remarking that a minimal deviation from condition~\eqref{eq:antonio2} already produces inconsistent behaviors, even extreme ones, showing that \eqref{eq:antonio2} is close to being necessary.




\bigskip

\section{Dependent processes for Gibbs--type priors}\label{sec:matteo}

In this section we briefly discuss possible extensions of the previous results to a dynamic setting. In particular, here we refer to time-indexed random objects, with some specification of the temporal transition mechanism, whose stationary, or at least marginal, states coincide in distribution with some random probability measure of Gibbs-type. In this respect, it is important to distinguish between two different research areas on time-dependent random probability measures, both related to Bayesian nonparametric priors.
The main difference between these two approaches, outlined below, lies in the fact that the former is mostly driven by inferential purposes, while the second is more concerned with the analytical properties of the constructed objects. If on one hand the first is closer to the interest of the Bayesian community, on the other it is our opinion that the two approaches have a strong potential of reciprocally benefitting from one another.

The first area, concerned with so-called dependent processes, is at present an extremely active front in Bayesian Nonparametrics.
Besides the pioneering contributions in \cite{CR78}, the modern approaches to the problem can be traced back to \cite{ME99}.
Generally speaking, the aim is to investigate generalizations of the Dirichlet process (or other random measures) to frameworks which allow for types of dependence less restrictive than exchangeability. These include for example dependence on time or, more generally, on covariates.  See, for example, \cite{BJQ12} for some up-to-date references. Most contributions in this direction exploit the representation (\ref{eq:discrete}) and dependence is quite easily induced via the weights and/or the atoms. Moreover, this allows to exploit simulation techniques such as the slice sampler (\cite{W07}, \cite{DWW99}) and the retrospective sampler \cite{PR08}. The combination of these two main factors leads then to efficient inferential procedures in such non exchangeable frameworks.

The second research area has its roots in Applied Probability and is concerned with stochastic population dynamics, but is also closely related to Bayesian nonparametric modeling. The main idea underlying the constructions in this framework is that of approximating the dynamics of a large population with a diffusion process, where the process dimension depends on the number of species the population is allowed to have. When the species can be of infinitely-many types, this gives rise to infinite-dimensional or measure-valued diffusions. In some cases the individual reproduction mechanisms yield populations whose frequencies have marginal or stationary states such as the one- and two-parameter Poisson-Dirichlet distribution (\cite{EK81},\cite{Pet(09)}), the Dirichlet process (\cite{EK86}), the normalized-inverse Gaussian distribution (\cite{RWF12}). From a Bayesian perspective these clearly represent dependent priors. At least in the authors' opinion, such an approach represents a highly promising research line for the definition of dependent processes, since the possibility of studying their analytical properties also yields a deeper understanding of their behavior. Other reasons of interest for the Bayesian community include the use of P\'olya urn schemes for constructing some of these dependent random probability measures (\cite{RW09},\cite{PR12}; see also \cite{CDD07}), and the investigation of the so-called $\sigma$-diversity processes (in the notation of Section \ref{sec:ramses}). These constitute a dynamic counterpart of \eqref{eq:kappan_asympt}, and make explicit the dynamics and distributional properties concerning the evolution of the clustering structure within the population, as a consequence of the specific modeling dynamics at hand. See \cite{RWF12} and \cite{R12}.

%
%

\bigskip

\section{Concluding remarks}\label{sec:conclusions}
An intense research activity, started after the introduction of the Dirichlet process, has produced a vast literature concerning classes of random probability measures whose laws can be used as nonparametric priors. In current research the choice among these classes is often dictated by taste (one's ``favorite prior''), mathematical tractability or a blend of the two. For instance, neutral to the right priors (\cite{dok(74)}) are typically used in survival analysis contexts since they are conjugate also w.r.t.~right censored observations. However, there is no conceptual reason to prefer a conjugate prior over a non--conjugate one and it all boils down to mathematical convenience since it allows to evaluate posterior inferences of interest. 
With Gibbs--type priors things go the opposite way: one makes a precise assumption on the learning mechanism according to which the prediction of a new value depends on the sample size $n$ and on the number of distinct values observed so far $K_n=k$ but not on their frequencies $n_1, \ldots, n_{k}$, and only afterwards investigates the implications of such an assumption. This is very much in the spirit of de~Finetti himself who constantly emphasizes in his works the importance of formulating assumptions on empirically 
``observable'' rather than on ``unobservable'' quantities.
In this respect Gibbs--type priors can be seen somehow as counterparts to characterizations of parametric families in terms of exchangeability and some other characteristic of the observables. Consider, for instance, Freedman's characterization \cite{Fre(62)} of exchangeable and rotational invariant sequences as mixtures of Gaussians: it is the request of rotational invariance on the observables that justifies the use of Gaussian distributions. In a nonparametric context, an analogous type of result (see \cite{Reg(78),Lo(91)}) legitimates the use of the Dirichlet process: by assuming exchangeability and a prediction rule given by a linear combination of the prior guess and the empirical measure one automatically obtains the Dirichlet process.

Turning back to the Gibbs--case, once the assumption on the learning mechanism is made, one realizes that the high degree of mathematical tractability is nothing but an implication and not a motivation. 
This then allows to work out a wealth of results concerning the behavior of Gibbs--type priors. Importantly, one is not anymore constrained to a logarithmic increase of $K_n$ as in the Dirichlet case and the whole spectrum going from a finite $K_n$ to an almost linearly increasing $K_n$ is available. This, in turn, produces a significantly more flexible prior on the number of components in mixture models. Furthermore, distributional properties and (often) closed form expressions for estimators of the quantities of statistical interest can be derived. An appealing feature is also represented by the fact that such quantities retain an intuitive flavor by directly relating to the key learning assumption. For instance, in the context of species sampling, one coherently has that $K_n$ is a sufficient statistic for predictions concerning ``new'' values. In contrast, if predictions are required for both ``new'' values and already observed values with frequency less than or equal to $\tau$, the sufficient statistics becomes $(K_n, M_{1,n}, \ldots, M_{\tau, n})$, which includes species with frequency not larger than $\tau$. Although more subtle, this is also in accordance with the key learning assumption and the implied reinforcement mechanism, described in the paper. Moreover, given the sound assumption on the learning scheme and its persuasive implications, it seems natural to use the Gibbs--framework also as basis for the definition of dependent processes.

Summing up, with this review we hope to have provided an affirmative and convincing answer to the question posed in the title of the paper. And we are confident that the future will see more statistical problems laid out in the well grounded general Gibbs--type framework. This would bring a solid foundation to the story and obviously would not prevent to use one's favorite Gibbs--type prior (e.g. the PY process) in the concrete application or even, if dropping the dependence on $K_n$ is legitimated by the problem at issue, returning to the ``safe'' Dirichlet world.


\bigskip

\section*{Acknowledgment}
The first three and last two authors are supported by the European Research Council (ERC) through StG "N-BNP" 306406. The fourth author is supported by CONACYT, project no. 131179.

\bigskip

\section*{Appendix}

\subsection*{Proof of Proposition \ref{prop:classification}} Recall that for any species sampling model the probability of generating a new value is of the form \eqref{eq:ss_new}. For it not to depend on $(n_1, \ldots, n_k)$ for any $n \ge 1$ and $k\leq n$, $p^{(n)}_{k}$ necessarily has to be of product form. \cite{Gne(05)} have shown that an EPPF associated to a infinite exchangeable random partition is of product form if and only if it is given by \eqref{eq:1} or, equivalently, if $\tilde p$ is of Gibbs--type. Therefore \eqref{eq:ss_new} depends only on $n$ and $k$ if and only if it is of Gibbs--type. This proves the categorization of species sampling models $\tilde p$ in classes (ii) and (iii).

We are now left with showing that the Dirichlet process is the only species sampling model for which \eqref{eq:ss_new} neither depends on the frequencies nor on $k$. Given the above, this amounts to showing that the subclass (i) of the family of Gibbs--type priors (ii) contains only the Dirichlet process.

First we show by a contradiction argument that for \eqref{eq:ss_new} not to depend on $k$ it must necessarily be $\sigma=0$. Then we conclude that the only Gibbs--type prior with $\sigma=0$ for which \eqref{eq:ss_new} does not depend on $k$ is the Dirichlet process.
Confining ourselves to the Gibbs--type case (since in all other cases \eqref{eq:ss_new} even depends on the frequencies), one has
     $$
    \mathbb{P}(X_{n+1}=\hbox{``new''} \: |\: X_1, \ldots, X_n)
    =\frac{V_{n+1,k+1}}{V_{n,k}} = 1 -(n-\sigma\, k) \frac{V_{n+1,k}}{V_{n,k}}
    $$
and we assume it does not depend on $k$. This amounts to requiring that
$
(n-\sigma\, k) V_{n+1,k} (V_{n,k})^{-1}
$
does not depend on $k$, namely
\begin{equation}\label{eq:old_no_dep}
\frac{V_{n+1,k}}{V_{n,k}}=\frac{c_n}{(n-\sigma\, k)},
\end{equation}
for some $c_n$ not depending on $k$ and, by using \eqref{eq:recursion},
\begin{equation}\label{eq:new_no_dep}
\frac{V_{n+1, k+1}}{V_{n,k}}=(1-c_n).
\end{equation}
The combination of \eqref{eq:old_no_dep} and \eqref{eq:new_no_dep} implies
\begin{equation}\label{eq:predict}
\mathbb{P}(X_{n+1} \in A \: |\: X_1, \ldots, X_n)= (1-c_n) P_0(A)+c_n \sum_{j=1}^k \frac{(n_j-\sigma)}{(n-\sigma \, k)} \delta_{X_j^*}(A).
\end{equation}
However, this is a prediction rule corresponding to an infinite exchangeable sequence if and only if $\sigma=0$. 
To see this, note that in view of \cite[Proposition 3.2]{For10} infinite exchangeability requires \eqref{eq:predict} to satisfy
\begin{equation}
  \label{eq:eugenio}
  \mathbb{P}(X_{n+1} \in A,\, X_{n+2}\in B\: |\: X_1, \ldots, X_n)=
  \mathbb{P}(X_{n+1} \in B,\, X_{n+2}\in A\: |\: X_1, \ldots, X_n)
\end{equation}
for any $n\ge 1$ and $A,B$ in $\mathscr{X}$. Consider, now, two sets $A$ and $B$ such that $A\cap B=\varnothing$, $A\cap\{X_1,\ldots,X_n\}=\varnothing$ and $B\cap \{X_1,\ldots,X_n\}\ne \varnothing$. Hence, the left--hand side of \eqref{eq:eugenio} coincides with
\[
c_n P_0(A)\left\{c_{n+1}P_0(B)+(1-c_{n+1})\frac{1}{n+1-(k+1)\sigma}\:
\sum_{j=1}^k(n_j-\sigma)\,\delta_{X_j^*}(B)\right\}
\]
whereas the right--hand side of \eqref{eq:eugenio} coincides with
\[
c_{n+1} P_0(A)\left\{c_{n}P_0(B)+(1-c_{n})\frac{1}{n-k\sigma}\:
\sum_{j=1}^k(n_j-\sigma)\,\delta_{X_j^*}(B)\right\}
\]
and the two are equal if and only if, for any $k=1,\ldots,n$, one has
\[
\frac{c_{n+1}(1-c_n)}{n-k\sigma}=\frac{c_n(1-c_{n+1})}{n+1-(k+1)\sigma}.
\]
Assuming $n\ge 2$, with $k=1$ the above condition becomes
\begin{equation}
  \label{eq:predictive1}
  \frac{c_{n+1}(1-c_n)}{n-\sigma}=\frac{c_n(1-c_{n+1})}{n+1-2\sigma}
\end{equation}
and with $k=2$, one has
\begin{equation}
  \label{eq:predictive2}
  \frac{c_{n+1}(1-c_n)}{n-2\sigma}=\frac{c_n(1-c_{n+1})}{n+1-3\sigma}.
\end{equation}
Taking the ratios of the terms in \eqref{eq:predictive1} and those in \eqref{eq:predictive2} yields
\[
\frac{n-2\sigma}{n-\sigma}=\frac{n+1-3\sigma}{n+1-2\sigma}
\]
and this holds true if and only if $\sigma^2=\sigma$. This therefore contradicts the assumption that $\mathbb{P}(X_{n+1}=\hbox{``new''} \: |\: X_1, \ldots, X_n)$ does not depend on $k$ for $\sigma \neq 0$. A different proof can also be derived by using the recursion \eqref{eq:recursion} iterated over two prediction steps.

Finally, recall that \cite[Theorem 13]{Gne(05)}, showed that Gibbs--type priors with $\sigma=0$ correspond to the Dirichlet process or the Dirichlet process mixture over its total mass parameter. 
On the other hand, when $\sigma=0$, \eqref{eq:predict} characterizes the Dirichlet process (see \cite{Lo(91)}, \cite{Reg(78)}). Hence, Dirichlet process mixture over the total mass cannot belong to class (i), i.e. $\mathbb{P}(X_{n+1}=\hbox{``new''} \: |\: X_1, \ldots, X_n)$ must depend also on $k$. The proof is, then, complete.
\hfill $\Box$

\end{document}